\documentclass{article}

\usepackage{amsmath}
\usepackage{amsthm}
\usepackage{amsxtra}

\usepackage{amsfonts,amssymb}
\usepackage{latexsym}

\newtheorem{Th}{Theorem}
\newtheorem{Prop}[Th]{Proposition}
\newtheorem{Lm}[Th]{Lemma}
\newtheorem{Co}[Th]{Corollary}

\theoremstyle{definition}
\newtheorem{Def}[Th]{Definition}

\newtheorem{Rem}{Remark}

\begin{title}
 {\bf Invariant states on the wreath product}
\end{title}
\author{ A.V. Dudko, N.I. Nessonov }

\date{}

\begin{document}
\maketitle
 \begin{abstract}
 Let $\mathfrak{S}_\infty$ be the infinity permutation group and
$\Gamma$ be a separable topological group.
 The wreath product $\Gamma\wr \mathfrak{S}_\infty$ is the
semidirect product $\Gamma^\infty_e \rtimes \mathfrak{S}_\infty$ for
the usual permutation action of $\mathfrak{S}_\infty$ on
$\Gamma^\infty_e=\left\{[\gamma_i]_{i=1}^\infty\,:\,\gamma_i\in
\Gamma,\textit{ only finitely many }\gamma_i\neq e\right\}$. In this
paper we obtain the full description of indecomposable states
$\varphi$ on the group $\Gamma\wr\mathfrak{S}_\infty,$ satisfying
the condition:
\begin{eqnarray*}
\varphi\left(sgs^{-1}\right)=
 \varphi\left( g\right)\text{ for each }g\in \Gamma\wr \mathfrak{S}_\infty,s\in\mathfrak{S}_\infty.
\end{eqnarray*}
 \end{abstract}

\subsubsection{Introduction}

\paragraph{ The  wreath product and $\mathfrak{S}_\infty$-central states.}
 Let $\mathbb{N}$ be the set of the natural numbers.
By definition,
 a bijection $s: \mathbb{N}\to \mathbb{N}$ is called
{\it finite} if the
set $\left\{ i\in\mathbb{N}|s(i) \neq i\right\}$ is finite.
 Define a group
$\mathfrak{S}_\infty$ as the group of all finite bijections
$\mathbb{N}\to\mathbb{N}$ and set
$\mathfrak{S}_n=\left\{s\in\mathfrak{S}_\infty|\; s(i)=i\;
 \text{ for each }
i>n\right\}$. Given a group $\Gamma$ identify
  element $\left(\gamma_1,\gamma_2,\ldots,\gamma_n\right)\in\Gamma^n$ with
 $\left(\gamma_1,\gamma_2,\ldots,\gamma_n,e\right)\in\Gamma^{n+1}$, where
$e$ is the identity element of $\Gamma$. The group $\Gamma^\infty_e$ is
defined as a inductive limit of sets
\begin{eqnarray}
\Gamma\mapsto\Gamma^2\mapsto\Gamma^3\mapsto\cdots\mapsto\Gamma^n\mapsto\cdots.
\end{eqnarray}
 The wreath product $\Gamma\wr \mathfrak{S}_\infty$ is the
semidirect product $\Gamma^\infty_e \rtimes \mathfrak{S}_\infty$
 for the usual
permutation action of $\mathfrak{S}_\infty$ on $\Gamma^\infty_e$.
 Using the
 imbeddings
 $\gamma\in \Gamma^\infty_e\to \left( \gamma,{\rm id}\right)
 \in \Gamma\wr \mathfrak{S}_\infty$,
 $s\in\mathfrak{S}_\infty\to\left( e^{(\infty)},s \right)\in
 \Gamma\wr \mathfrak{S}_\infty$, where $e^{(\infty)}=
 \left( e,e,\ldots,e,\ldots\right)$ and ${\rm id}$ is the
 identical bijection,
 we  identify $\Gamma^\infty_e$ and $\mathfrak{S}_\infty$
 with the corresponding
subgroups of $\Gamma\wr \mathfrak{S}_\infty$.
 Therefore, each element $g$ of $\Gamma\wr \mathfrak{S}_\infty$ is of the
form $g=s\gamma$, with $\gamma =\left(\gamma_1,
\gamma_2,\ldots\right)\in\Gamma^\infty_e$ and $s\in\mathfrak{S}_\infty$.
Furthermore, it is assumed that $s\left(\gamma_1, \gamma_2,\ldots
\right)s^{-1}=\left(\gamma_{s^{-1}(1)}, \gamma_{s^{-1}(2)},\ldots \right)$.

If $\Gamma$ is a
 topological group, then we
will equip $\Gamma^n$ with the natural
 product-topology.  Furthermore, we will always consider
$ \Gamma^\infty_e $ as a topological group with
 the inductive limit topology.
 The group $\Gamma\wr \mathfrak{S}_\infty$
 is isomorphic to  $ \Gamma^\infty_e \times
\mathfrak{S}_\infty $, as a set. Therefore, we will
 equip the group
 $\Gamma\wr \mathfrak{S}_\infty$ with the product-topology,
 considering
 $\mathfrak{S}_\infty$ as a discrete topological space. From now on we
 assume that $\Gamma$ is a separable topological group.

  \paragraph{ The basic definitions. }
  Let $\mathcal{H}$ be a Hilbert space, let
 $\mathcal{B}\left(\mathcal{H}\right)$ be the set of all
 bounded operators in $\mathcal{H}$ and let
 $\mathcal{I}_{\mathcal{H}}$ be the identity operator
 in $\mathcal{H}$.  We denote by $\mathcal{U}\left(
\mathcal{H} \right)$ the unitary subgroup in
 $\mathcal{B}\left(\mathcal{H}\right)$.
 By a unitary representation of the topological group
 $G$ we will always
mean a {\it continuous} homomorphism of $G$ into
$\mathcal{U}\left(\mathcal{H} \right)$, where
 $\mathcal{U}\left(\mathcal{H} \right)$
 is equipped with the strong operator topology. For unitary
 representation $\pi$ of the group $G$ we denote $\mathcal{M}_\pi$
 the $W^*-$algebra $\pi(G)^{\prime\prime}$, which is
generated by the operators $\pi(g)\; \left( g\in G\right)$.

 \begin{Def}\label{indecomposable}
 An unitary representation
 $\pi:G\to\mathcal{U}\left(\mathcal{H} \right)$ of
the group $G$ is called a factor-representation if
 $\mathcal{M}_\pi$
 is a factor. A positive definite function $\varphi$ on  group $G$
 is called an indecomposable, if the corresponding GNS-representation
  is a factor-representation.
 \end{Def}
Further, an element $\Gamma\wr\mathfrak{S}_\infty$
 can always be written
as the product of an element from $\mathfrak{S}_\infty$
 and an element
 from $\Gamma^\infty_e $. The commutation rule between
these two kinds of elements is
 \begin{equation}\label{product}
s \gamma=s \left(\gamma_1,\gamma_2,\ldots \right)=
 \left(
\gamma_{s^{-1}(1)},\gamma_{s^{-1}(2)},\ldots \right)  s,
\end{equation}
where $s\in\mathfrak{S}_\infty$,$\gamma=\left(\gamma_1, \gamma_2,\ldots
\right)\in\Gamma^\infty_e $. Let $\mathbb{N}\diagup s$
 be the set of orbits
of $s$ on the set $\mathbb{N}$. Note that for $p\in \mathbb{N}\diagup s$
 permutation $s_p$, which is
 defined by the formula
 \begin{equation}\label{sp}
 s_p(k)=
\left\{
 \begin{array}{rl}
 s(k)&\text{ if } k\in p\\
 k&\textit{otherwise}
 \end{array}\right.,
  \end{equation}
 is a cycle of the order $|p|$, where $|p|$ denotes the
 cardinality of $p$.
For $\gamma=\left(\gamma_1,\gamma_2,\ldots \right)
 \in\Gamma^\infty_e$
we define the element $\gamma(p)=\left(\gamma_1(p), \gamma_2(p),\ldots
\right) \in\Gamma^\infty_e $ as follows
 \begin{equation}\label{color1}
 \gamma_k(p)= \left\{
 \begin{array}{rl}
 \gamma_k&\text{ if } k\in p\\
 e&\textit{otherwise}.
 \end{array}\right.
  \end{equation}
 Thus, using ({\ref{product}}), we have
 \begin{equation}\label{decompositiontocycles}
s \gamma = \prod\limits_{p\in \mathbb{N}\diagup s}
 s_p \gamma(p).
\end{equation}
Element $ s_p \gamma(p)$ is called the {\it generalized cycle} of
$s\gamma$.

 Denote by $(n\,\;k)\in\mathfrak{S}_\infty$ the transposition of
numbers $k$ and $n.$ Following Olshanski (see \cite{O2}) we
introduce permutations
$\omega_n=\omega^{(0)}_n\in\mathfrak{S}_\infty$ by the next formula:
\begin{eqnarray}\label{omega}
\omega_n(i)=\left\{\begin{array}{ll}i,&\textit{ if }\;2n<i,\\
i+n,&\textit{ if }\;i\leqslant n,\\
i-n,&\textit{ if }\;n<i\leqslant 2n.\end{array}\right.
\end{eqnarray}

For the element $g=s  \gamma$ we call $\textit{support}$ of $g$ the
set ${\rm supp}(g)=\left\{i:s(i)\neq i\text{ or }\gamma_i\neq
e\right\}$. Note that ${\rm supp}(g)$ is always finite subset of
$\mathbb{N}$. If ${\rm supp}(g_1)\cap {\rm supp}(g_2)=\emptyset$
then elements $g_1$ and $g_2$ commute.

\begin{Def}\label{def central}
Let $G$ be a group and let $H$ be a subgroup of $G$. A positive definite
function $\varphi$ on
 $G$ is called $H$-central  if  $\varphi(gh)=\varphi(hg)$ for all
 $h\in H$ and $g\in G$. We say that $\varphi$ is a {\it state} on $G$, if
$\varphi(e)=1$, where $e$ is the identical element of $G$. A state $\varphi
$ is called {\it indecomposable}, if the corresponding GNS-representation
$\pi_\varphi$ is a factor representation.
\end{Def}
Let $\mathcal{M}_*$ denotes the space of all $\sigma$-weakly continuous
functional on $w^*$-algebra $\mathcal{M}$.

Now we fix a $\mathfrak{S}_\infty$-central state $\varphi$ on $\Gamma\wr
\mathfrak{S}_\infty$, and denote by $\pi_\varphi$ the corresponding
GNS-representations.
\begin{Th}\label{Th3}
Let $\pi_\varphi\left( \Gamma\wr \mathfrak{S}_\infty\right)^{\prime\prime}$
be a $w^*$-algebra generated by operators $\pi_\varphi\left(\Gamma\wr
\mathfrak{S}_\infty \right)$ and let $\mathcal{C}\left(
\pi_\varphi\left(\Gamma\wr \mathfrak{S}_\infty \right)\right)$ be the center of $\pi_\varphi\left( \Gamma\wr \mathfrak{S}_\infty\right)^{\prime\prime}$.
Suppose that the positive functionals $\varphi_1$ and $\varphi_2$ from
 $\pi_\varphi\left( \Gamma\wr \mathfrak{S}_\infty\right)^{\prime\prime}_*$
 satisfy the next conditions:
 \begin{itemize}
\item {\rm i}) $\varphi_k\left(\pi_\varphi(s) a \right) =
 \varphi_k\left( a \pi_\varphi(s)\right)$ for all $s\in\mathfrak{S}_\infty$
and $a\in\pi_\varphi\left( \Gamma\wr
\mathfrak{S}_\infty\right)^{\prime\prime}$ $\left(k=1,2 \right)$;
 \item
 {\rm ii}) $\varphi_1\left(\mathfrak{c} \right) = \varphi_2\left(\mathfrak{c}
\right) $ for all $\mathfrak{c}\in \mathcal{C}\left(
\pi_\varphi\left(\Gamma\wr \mathfrak{S}_\infty \right)\right)$.
\end{itemize}
 Then $\varphi_1\left(\mathfrak{a} \right) = \varphi_2\left(\mathfrak{a}
\right) $ for all $\mathfrak{a}\in \pi_\varphi\left(\Gamma\wr
\mathfrak{S}_\infty \right)$.
\end{Th}
Recall that representations $\pi_1$ and $\pi_2$ of the group $G$ are called
quasiequiva\-lent if there exists isomorphism $\theta:\pi_1\left(G
\right)^{\prime\prime}\mapsto\pi_2\left(G \right)^{\prime\prime}$ with the
property
\begin{eqnarray}
\theta\left(\pi_1\left(g \right)\right)=\pi_2\left(g \right) \text{ for all
} g\in G.
\end{eqnarray}
The following corollary is immediate consequence of the above theorem.
\begin{Co}\label{Co4}
If $\varphi_1$ and $\varphi_2$ are indecomposable
$\mathfrak{S}_\infty$-central states on $\Gamma\wr \mathfrak{S}_\infty$
such that the corresponding GNS-representations $\pi_{\varphi_1}$ and
$\pi_{\varphi_2}$ are quasiequivalent, then $\varphi_1=\varphi_2$.
\end{Co}
\paragraph{The natural examples.}\label{parnatexmmpl}
For any state $\varphi $ on $\Gamma$ define two
$\mathfrak{S}_\infty$-central states $\varphi_{sp} $ and
$\varphi_{reg}$ on $\Gamma\wr \mathfrak{S}_\infty$ as follows
\begin{eqnarray}
\varphi_{sp} \left(s\gamma \right)=\prod\varphi \left(\gamma_k \right)
\text{ for all } \gamma=\left(\gamma_1,\gamma_2,\ldots
\right)\in\Gamma^\infty_e \text{ and } s\in\mathfrak{S}_\infty;\label{sp0}\\
 \varphi_{reg}\left(s\gamma \right)=\left\{
 \begin{array}{rl}
 \prod\varphi \left(\gamma_k \right)&\text{ if }  s=e\\
 0&\text{ if } s\neq e.
 \end{array}\right.
\end{eqnarray}
We have the following result:
\begin{Prop}
For GNS-representations $\pi_{\varphi_{sp} }$ and $\pi_{\varphi_{reg} }$
the next properties hold:
\begin{itemize}
  \item {\rm(i)} If $\pi_{\varphi_{sp} }$ acts in Hilbert space
      $\mathcal{H}_{\varphi_{sp} }$, and
      $\mathcal{H}_{\varphi_{sp}}^\mathfrak{S}=\left\{\eta\in\mathcal{H}_{\varphi_{sp}
      }: \pi_{sp}(s)\eta=\eta \text{ for all } s\in \mathfrak{S}_\infty
      \right\}$, then ${\rm dim}
      \mathcal{H}_{\varphi_{sp}}^\mathfrak{S}=1 $. In particular,
      $\pi_{\varphi_{sp} }$ is irreducible.
\item {\rm(ii)}  $\pi_{\varphi_{reg} }$ is a factor representation.
    \item {\rm(iii)}  $w^*$-algebra $\pi_{\varphi_{reg} }
    \left(\Gamma\wr \mathfrak{S}_\infty\right)^{\prime\prime}$ is a
    factor of the  type ${\rm II}$ or ${\rm III}$.
\end{itemize}
\end{Prop}
\begin{proof}
 Let $\xi_{\varphi_{sp}}$ $\left(\xi_{\varphi_{reg}} \right)$ be the cyclic
 vector for representation $\pi_{sp}$  $\left(\pi_{reg}\right)$ with the
property
\begin{eqnarray*}
  &\varphi_{sp}\left(g \right)=\left(\pi_{sp}(g)
 \xi_{\varphi_{sp}},\xi_{\varphi_{sp}} \right)&\;\;
 \bigg(\varphi_{reg}\left(g \right)=\left(\pi_{reg}(g)
 \xi_{\varphi_{reg}},\xi_{\varphi_{reg}} \right)\bigg)\\
&\text{ for all }
 g\in \Gamma\wr
\mathfrak{S}_\infty.&
 \end{eqnarray*}
Set $\Gamma^{n \infty}_e= \left\{\gamma=\left(\gamma_1,\gamma_2,\ldots
\right)\in\Gamma^\infty_e\big| \gamma_k=e\text{ for all } k\leq
n\right\}$,\\
 $\mathfrak{S}_{n \infty}= \left\{s\in\mathfrak{S}_{\infty}\big|
 s(k)=k \text{ for all } k\leq n\right\}$. Denote by $\Gamma\wr
\mathfrak{S}_{n\infty}$ the subgroup of
 $\Gamma\wr\mathfrak{S}_\infty$ generated by $\Gamma^{n \infty}_e$ and
 $\mathfrak{S}_{n \infty}$.

 To the proof point {\rm(i)}, first  we note that, by definition GNS-construction,
 $\xi_{\varphi_{sp}}$ lies in $\mathcal{H}_{\varphi_{sp}}^\mathfrak{S}$.
Further we will use the important mixing-property. Namely, denote by
$\omega_n$ a bijection which acts as follows
\begin{eqnarray}\label{omegan}
\omega_n(i)=\left\{\begin{array}{ll}i,&\textit{ if }\;2n<i,\\
i+n,&\textit{ if }\;i\leqslant n,\\
i-n,&\textit{ if }\;n<i\leqslant 2n.\end{array}\right.
\end{eqnarray}
Then for any $\eta \in \mathcal{H}_{\varphi_{sp}}^\mathfrak{S}$, using
(\ref{sp0}), we obtain
\begin{eqnarray}
\lim\limits_{n\to\infty}\left(\pi_{sp}\left(\omega_n \right)\eta,\eta
\right) =\left(\xi_{\varphi_{sp}},\eta \right)\left(\eta,\xi_{\varphi_{sp}}
\right).
\end{eqnarray}
This implies {\rm (i)}.

 A property {\rm (ii)} follows from Proposition
\ref{multiplicativity} (below). Nevertheless, using the explicit
realizations of  $\pi_{\varphi_{reg} }$, we give another proof. We begin
with the GNS-representation $T$ of $\Gamma$ which acts in Hilbert space
$\mathcal{H}_T$ with cyclic vector $\xi_\varphi$: $\varphi \left(\gamma
\right)=\left(T(\gamma )\xi_\varphi,\xi_\varphi\right)$ for all
$\Gamma\in\gamma$. Further, using embedding $\mathcal{H}_T^{\otimes n}\ni
\eta\mapsto \eta\otimes\xi_\varphi\in\mathcal{H}_T^{\otimes n+1}$, define
Hilbert space $\mathcal{H}_T^{\otimes \infty}$ and corresponding
representation $T^{\otimes\infty}$ of $\Gamma^\infty_e$:
\begin{eqnarray*}
T^{\otimes\infty}(\gamma) \left(\xi_1\otimes\xi_2\otimes\ldots \right)=
T\left(\gamma_1 \right)\xi_1\otimes T\left(\gamma_2
\right)\xi_2\otimes\ldots, \text{ where } \gamma
=\left(\gamma_1,\gamma_2,\ldots \right).
\end{eqnarray*}
The action $U$ of $\mathfrak{S}_\infty$ on $\mathcal{H}_T^{\otimes \infty}$
is given by the formula
\begin{eqnarray*}
U(s)\left(\xi_1\otimes\xi_2\otimes\ldots\otimes\xi_k \otimes\ldots\right)=
\xi_{s^{-1}(1)}\otimes\xi_{s^{-1}(2)}\otimes\ldots\otimes\xi_{s^{-1}(k)}
\otimes\ldots
\end{eqnarray*}
 Now we define operator $\Pi(g)$
$\left(g\in\Gamma\wr\mathfrak{S}_\infty\right)$ in $l^2\left(
\mathfrak{S}_\infty,\mathcal{H}_T^{\otimes \infty}\right)$ as follows
\begin{eqnarray*}
\left(\Pi(\gamma)\eta\right)(s)=U(s)T^{\otimes\infty}(\gamma)U^*(s)\eta(s)\;\;
 \left(\gamma \in \Gamma^\infty_e, \eta\in l^2\left(
\mathfrak{S}_\infty,\mathcal{H}_T^{\otimes \infty}\right) \right);\\
\left(\Pi(t)\eta\right)(s)=\eta(st) \;\;\left(t\in\mathfrak{S}_\infty \right).
\end{eqnarray*}
Since for any $s\in\mathfrak{S}_\infty$ and $g=\left(\gamma_1,
\gamma_2,\ldots \right)\in\Gamma_e^\infty$ $s\left(\gamma_1,
\gamma_2,\ldots \right)s^{-1}=\left(\gamma_{s^{-1}(1)},
\gamma_{s^{-1}(2)},\ldots \right)$, $\Pi$ extends by multiplicativity to
the representation of $\Gamma\wr\mathfrak{S}_\infty$.

 If  ${\xi}_\varphi^{\otimes\infty}
=\xi_\varphi\otimes\xi_\varphi\otimes\ldots\in\mathcal{H}_T^{\otimes
\infty}$ and $\widehat{\xi}_\varphi (g)=\left\{\begin{array}{ll}
{\xi}_\varphi^{\otimes\infty},
&\textit{ if }\;g=e,\\
0,&\textit{ if }\;g\neq e\end{array}\right. $ then we have
\begin{eqnarray}
\varphi_{reg} \left(s\gamma\right)=\left(\Pi(s\gamma
)\widehat{\xi}_\varphi,
 \widehat{\xi}_\varphi\right) \;\; \left(s\in\mathfrak{S}_\infty,\gamma
\in\Gamma^\infty_e\right).
\end{eqnarray}
Therefore, without loss generality we can assume that $\pi_{reg}=\Pi$.

Let $\Pi^\prime$ denote the representation of $\mathfrak{S}_\infty$ which
acts on $l^2\left( \mathfrak{S}_\infty,\mathcal{H}_T^{\otimes
\infty}\right)$ by
\begin{eqnarray}
\left(\Pi^\prime(t)\eta\right)(s)=U(t)\eta(t^{-1}s).
\end{eqnarray}
Obvious, $\Pi^\prime\left(\mathfrak{S}_\infty\right)$ is contained in
commutant $\Pi\left(\Gamma\wr \mathfrak{S}_{\infty}\right)^\prime$ of
 $\Pi\left(\Gamma\wr \mathfrak{S}_{\infty}\right)$.

Let us prove that center $\mathcal{C}=\Pi\left(\Gamma\wr
\mathfrak{S}_{\infty}\right)^{\prime\prime}\cap\Pi\left(\Gamma\wr
\mathfrak{S}_{\infty}\right)^\prime$ of $\Pi\left(\Gamma\wr
\mathfrak{S}_{\infty}\right)^{\prime\prime}$ is trivial.

Our proof starts with the observation that
\begin{eqnarray}
\Pi(g)\Pi^\prime(g) \widehat{\xi}_\varphi=\widehat{\xi}_\varphi\text{ for
all } g\in\mathfrak{S}_\infty.
\end{eqnarray}
Hence for $\mathfrak{c}\in\mathcal{C}$ we have
\begin{eqnarray}\label{fixed}
\Pi(g)\Pi^\prime(g)\mathfrak{c}
\widehat{\xi}_\varphi=\mathfrak{c}\widehat{\xi}_\varphi\text{ for all }
g\in\mathfrak{S}_\infty.
\end{eqnarray}
In particular, this gives
\begin{eqnarray}
\left\| \mathfrak{c}\widehat{\xi}_\varphi(s)\right\|=
 \left\| \mathfrak{c}\widehat{\xi}_\varphi\left(gsg^{-1} \right)\right\|
\text{ for all } g,s\in\mathfrak{S}_\infty.
\end{eqnarray}
Since every conjugacy class $C(s)= \left\{gsg^{-1}:g\in
\mathfrak{S}_\infty\right\}$ is infinite except $s=e$, we have
\begin{eqnarray}
\mathfrak{c}\widehat{\xi}_\varphi(s)=0 \text{ for all } s\neq e.
\end{eqnarray}
It follows from (\ref{fixed}) that
\begin{eqnarray}
U(s)\left(\mathfrak{c}\widehat{\xi}_\varphi(e) \right)=
\mathfrak{c}\widehat{\xi}_\varphi(e) \text{ for all }
s\in\mathfrak{S}_\infty.
\end{eqnarray}
As in the proof of the point {\rm(i)}, this gives that
$\mathfrak{c}\widehat{\xi}_\varphi(e)=\alpha {\xi}_\varphi^{\otimes\infty}$
 $\left(\alpha \in\mathbb{C} \right)$. Since $\widehat{\xi}_\varphi$ is
cyclic, we have $\mathfrak{c}=\alpha I$. Therefore, $w^*$-algebra
 $\Pi\left(\Gamma\wr
\mathfrak{S}_{\infty}\right)^{\prime\prime}$ is a factor.

{\rm (iii)} We begin by recalling the notion of a {\it central sequence} in
a factor $\mathcal{M}$. A bounded sequence $ \left\{a_n
\right\}\subset\mathcal{M}$ is called {\it central} if
\begin{eqnarray*}
 s-\lim\limits_{n\to\infty}\left(a_n m -ma_n\right)=0\text{ and }
s-\lim\limits_{n\to\infty}\left(a_n^* m -ma_n^*\right)=0 \text{ for all
}m\in\mathcal{M}.
\end{eqnarray*}
A {\it central sequence} is called {\it trivial} if there exists sequence $
\left\{c_n \right\}\subset\mathbb{C}$ such that
\begin{eqnarray*}
 s-\lim\limits_{n\to\infty}\left(a_n  -c_nI\right)=0\text{ and }
s-\lim\limits_{n\to\infty}\left(a_n^* -\overline{c}_n I\right)=0.
\end{eqnarray*}
Let $s_k$ be the transposition interchanging $k$ and $k+1$. We claim that $
\left\{\pi_{reg}\left(s_n \right)\right\}$ is non trivial cental sequence.
Indeed, since $\varphi_{reg}$ is a $\mathfrak{S}_\infty$-central state, we
have
 \begin{eqnarray*}
\lim\limits_{n\to\infty}\left(m\pi_{reg}\left(s_n
\right)-\pi_{reg}\left(s_n \right)m\right){\xi_{\varphi_{reg} }}=0 \text{
for all } m\in
 \Pi\left(\Gamma\wr
\mathfrak{S}_{\infty}\right)^{\prime\prime}.
\end{eqnarray*}
It follows that
\begin{eqnarray*}
\lim\limits_{n\to\infty}\left(m\pi_{reg}\left(s_n
\right)-\pi_{reg}\left(s_n \right)m\right)x{\xi_{\varphi_{reg} }}=0 \text{
for all } m,x\in
 \Pi\left(\Gamma\wr
\mathfrak{S}_{\infty}\right)^{\prime\prime}.
\end{eqnarray*}
Since $\xi_{\varphi_{reg} } $ is cyclic and $\varphi_{reg}\left(s_n
\right)=0$, then $ \left\{\pi_{reg}\left(s_n \right)\right\}$ is non
trivial central sequence.

It remains to prove that each cental sequence in factor $\mathcal{M}$ of
type ${\rm I}$ is trivial.
 Suppose that $\mathcal{M}$ is a factor of type ${\rm I}$.
 Let $\left\{\mathfrak{e}_{kl}: k,l\in\mathbb{N} \right\}$
 be a matrix unit in $\mathcal{M}$. This means that the next relations hold
\begin{eqnarray}\label{munit}
\mathfrak{e}_{kl}^*=\mathfrak{e}_{lk},\;
\mathfrak{e}_{kl}\mathfrak{e}_{pq}=\delta_{lp}\mathfrak{e}_{kq},\;
\sum\limits_{k\in\mathbb{N}} \mathfrak{e}_{kk}=I.
\end{eqnarray}
Let $ \left\{a_n=\sum\limits_{k,l}c_{kl}(n)
\mathfrak{e}_{kl}:c_{kl}(n)\in\mathbb{C}\right\}$ be a cental sequence in
$\mathcal{M}$. Set
$\mathfrak{C}_{pq}(n)=a_n\mathfrak{e}_{pq}-\mathfrak{e}_{pq}a_n$. An easy
computation shows that
\begin{eqnarray*}
\mathfrak{e}_{qq}\left(\mathfrak{C}_{pq}(n)\right)^*\mathfrak{C}_{pq}(n)
\mathfrak{e}_{qq}=
\left[\left|c_{pp}(n)-c_{qq}(n) \right|^2-\left|c_{pp}(n)
\right|^2+\sum\limits_k\left|c_{kp}(n) \right|^2\right]\mathfrak{e}_{qq},\\
\mathfrak{e}_{pp}\mathfrak{C}_{pq}(n)
\left(\mathfrak{C}_{pq}(n)\right)^*\mathfrak{e}_{pp}=
 \left[\left|c_{pp}(n)-c_{qq}(n) \right|^2-\left|c_{qq}(n)
\right|^2+\sum\limits_k\left|c_{qk}(n) \right|^2\right]\mathfrak{e}_{pp}.
\end{eqnarray*}
Using the fact that $\left\{a_n\right\}$ is a central sequence, we deduce
from this that
\begin{eqnarray*}
\lim\limits_{n\to\infty}\sum\limits_{k:k\neq q}\left|c_{qk}(n)
\right|^2=0,\; \lim\limits_{n\to\infty}\sum\limits_{k:k\neq
q}\left|c_{kq}(n) \right|^2=0,\\
 \lim\limits_{n\to\infty}\left|c_{11}(n)
-c_{qq}(n)\right|^2=0 \text{ for all } q.
\end{eqnarray*}
This means that $s-\lim\limits_{n\to\infty}\left( a_n-c_{11}(n)I\right)=0$
and $s-\lim\limits_{n\to\infty}\left(
a_n^*-\overline{c_{11}(n)}I\right)=0$. Thus $\left\{a_n\right\}$ is
trivial.
\end{proof}

The goal of this paper is to give the full description of indecomposable
$\mathfrak{S}_\infty$-central states on $\Gamma\wr\mathfrak{S}_\infty$ (see definition \ref{def central}).
The character theory of infinite wreath product in the case  of finite $\Gamma$ is
developed by R. Boyer \cite{Boy}. In this case $\Gamma\wr\mathfrak{S}_\infty$
is inductive limit of finite groups, their finite characters can be obtained
 as limits of normalized characters of prelimit finite groups,
 and Boyer's method is a direct generalization of Vershik's-Kerov's asymptotic
 approach \cite{VK0}. The characters of $\Gamma\wr\mathfrak{S}_\infty$ for
 general separable group $\Gamma$ were found by authors in \cite{DN},
 \cite{DN1}. Our method has been based on the ideas of Okounkov, which
 he has developed for the proof of Thoma's theorem \cite{Thoma},
 \cite{Ok1}, \cite{Ok2}.

 A finite character is a  $\Gamma\wr\mathfrak{S}_\infty$-central
  positive definite function on
   $\Gamma\wr\mathfrak{S}_\infty$. In this paper we study the more general
   class of the $\mathfrak{S}_\infty$-central states on
   $\Gamma\wr\mathfrak{S}_\infty$. Our results provide a complete classification
   such indecomposable states. The set of all indecomposable  $\mathfrak{S}_\infty$-central states
   have very important property. Namely, if for for two indecomposable
   $\mathfrak{S}_\infty$-central states $\varphi_1$  and $\varphi_2$ the corresponding
    GNS-representations $\pi_{\varphi_1}$ and
 $\pi_{\varphi_2}$ are quasiequivalent, then $\varphi_1=\varphi_2$
 (theorem \ref{Th3}, corollary \ref{Co4}).

 The papers is organized as follows. Below we give a brief description  of
 the general properties of the $\mathfrak{S}_\infty$-central states.
 The key results are
 lemma \ref{weak-lim} and proposition \ref{multiplicativity}.  Here we also
 recall the classification of the traces (central states) on
 $\Gamma\wr\mathfrak{S}_\infty$ (theorem \ref{mainth}). In section
 \ref{ex of repr} we present the full collection of  factor-representations,
  which  define the $\mathfrak{S}_\infty$-central states
   (proposition \ref{Prop11a}).
   Each such state is parametrized by pair $\left(A,\rho\right)$, where
   $A$ is self-adjoint operator, $\rho$ is the unitary representation
    of $\Gamma $ (paragraph \ref{paragraph2.1}). In proposition \ref{Prop12} we prove that the unitary
    equivalence of pairs
    $\left(A_1,\rho_1 \right)$ and $\left(A_2,\rho_2 \right)$ is
    equivalent  to the equality of the corresponding $\mathfrak{S}_\infty$-central
    states. In section \ref{KMSsec} we discuss about physical  KMS-condition (see \cite{Tak})
    for these states (theorem \ref{theorem15}). In section \ref{Themainresult} we prove
     the classification theorem \ref{main}.

\paragraph{ The multiplicativity.}
Let $\varphi$ be an indecomposable $\mathfrak{S}_\infty$-central
state on the group $\Gamma\wr \mathfrak{S}_\infty$. Then it defines
according to GNS-construction a factor-representation $\pi_\varphi$
of the group $\Gamma\wr \mathfrak{S}_\infty$ with cyclic vector
$\xi_\varphi$ such that
$\pi_\varphi(g)=\left(\pi_\varphi(g)\xi_\varphi,\xi_\varphi\right)$
for each $g\in \Gamma\wr \mathfrak{S}_\infty$. The next lemma shows,
that different indecomposable $\mathfrak{S}_\infty$-central states
define representations which are not quasiequivalent. Let $ w-lim$
stand for the limit in the weak operator topology.
\begin{Lm}\label{weak-lim}
Let $\varphi$ be an indecomposable $\mathfrak{S}_\infty$-central
state on the group $\Gamma\wr \mathfrak{S}_\infty$. Than for each
$g\in \Gamma\wr \mathfrak{S}_\infty$ there exists
$w-\lim\limits_{n\rightarrow\infty}\pi_\varphi\left(\omega_n
g\omega_n\right)$ and the next equality holds:
\begin{eqnarray}\label{weak-lim equality}
w-\lim\limits_{n\rightarrow\infty}\pi_\varphi\left(\omega_n
g\omega_n\right)=\varphi(g)  I.
\end{eqnarray}
\end{Lm}
\begin{proof}
   Let $h_1,h_2\in \Gamma\wr \mathfrak{S}_\infty$. Fix $k$ such that
  \begin{eqnarray}\label{supp g}
  {\rm supp}(h_1),{\rm supp}(h_2),{\rm supp}(g)\subset\{1,2,\ldots,k\}.\end{eqnarray} For each
  $n\in \mathbb{N}$ there exists elements
  $g_{(n,k)},h_{(n,k)}\in \mathfrak{S}_\infty$ such that
  \begin{eqnarray}\label{supp gnk}
  {\rm supp}(g_{(n,k)}),{\rm supp}(h_{(n,k)})\subset \{k+1,k+2,\ldots\}\end{eqnarray}
  and
  $\omega_{n+k}=g_{(n,k)}\omega_kh_{(n,k)}$ (see (\ref{omega})). Permutations $g_{(n,k)},h_{(n,k)}$ can
  be defined as follows:
  \begin{eqnarray*}
g_{(n,k)}(i)=\left\{\begin{array}{ll}
i,&\textit{ if }\;i\leqslant k\textit{ or }2k+2n<i,\\
i+n,&\textit{ if }\;k<i\leqslant 2k+n,\\
i-k-n,&\textit{ if }\;2k+n<i\leqslant 2k+2n.\end{array}\right.\\
h_{(n,k)}(i)=\left\{\begin{array}{ll}
i,&\textit{ if }\;i\leqslant k\textit{ or }2k+n<i,\\
i+k,&\textit{ if }\;k<i\leqslant k+n,\\
i-n,&\textit{ if }\;k+n<i\leqslant 2k+n.\end{array}\right.
\end{eqnarray*}
By (\ref{supp g}) and (\ref{supp gnk}), the elements $g_{(n,k)}$ and
$h_{(n,k)}$ commutes with the elements $h_1,h_2$ and $g$. Therefore
\begin{eqnarray}\label{h(n,k)}\begin{split}
&h_2^{-1}\omega_{n+k}
 g\omega_{n+k}h_1=h_2^{-1}
 \left(g_{(n,k)}\omega_kh_{(n,k)}\right)^{-1}
   g  g_{(n,k)}\omega_kh_{(n,k)}
   h_1\\&=h_{(n,k)}^{-1}h_2^{-1}\omega_k
 g\omega_kh_1h_{(n,k)}.
 \end{split}
\end{eqnarray}As $\varphi$ is $\mathfrak{S}_\infty$-central, one has:
\begin{eqnarray}\label{stabilisation}
\begin{split}
&\left(\pi_\varphi\left(\omega_{n+k}
 g\omega_{n+k}\right)\pi_\varphi(h_1)\xi_\varphi,\pi_\varphi(h_2)\xi_\varphi\right)
 =\varphi\left(h_2^{-1}\omega_{n+k}
 g\omega_{n+k}h_1\right)=\\
 &\varphi\left(h_2^{-1}\omega_k
 g\omega_kh_1\right)=
 \left(\pi_\varphi\left(\omega_k
 g\omega_k\right)\pi_\varphi(h_1)\xi_\varphi,\pi_\varphi(h_2)\xi_\varphi\right).
 \end{split}
\end{eqnarray}
 As $\xi_\varphi$ is cyclic, by
 $(\ref{stabilisation})$,
 there exists the limit
  \begin{eqnarray*}w-\lim\limits_{n\rightarrow\infty}\pi_\varphi(\omega_n
g\omega_n).\end{eqnarray*} For each $h\in \Gamma\wr
\mathfrak{S}_\infty$ for large enough $n$ one has ${\rm
supp}(\omega_ng\omega_n)\cap {\rm supp}(h)=\emptyset$. Therefore
$\pi_\varphi(\omega_ng\omega_n)\pi_\varphi(h)=\pi_\varphi(h)\pi_\varphi(\omega_ng\omega_n)$.
This involves that the weak limit
$w-\lim\limits_{n\rightarrow\infty}\pi_\varphi(\omega_n g\omega_n)$
lies in the center of the algebra $M_{\pi_\varphi}$, generated by
operators of
 the representation $\pi_\varphi$. Thus $\lim\limits_{n\rightarrow\infty}\pi_\varphi(\omega_n g\omega_n)$ is scalar.
 By
 $\mathfrak{S}_\infty$-centrality of $\varphi,$
 \begin{eqnarray*}\left(w-\lim\limits_{n\rightarrow\infty}\pi_\varphi(\omega_n
g\omega_n)\xi_\varphi,\xi_\varphi\right)=\lim\limits_{n\rightarrow\infty}\varphi(\omega_n
g\omega_n)=\varphi(g),\end{eqnarray*} which finishes the proof.
 \end{proof}
 The following claim gives a useful characterization
 of the class of the indecomposable $\mathfrak{S}_\infty$-central states:
 \begin{Prop}\label{multiplicativity}
The following conditions for $\mathfrak{S}_\infty$-central state
$\varphi$ on the group
  $ \Gamma\wr\mathfrak{S}_\infty$ are equivalent:
 \begin{itemize}
\item [\it (a)] $\varphi$ is indecomposable;
 \item [\it (b)] $\varphi(gg')=\varphi(g)\varphi(g')$ for each
 $g,g'\in\Gamma\wr\mathfrak{S}_\infty$ with
 ${\rm supp}(g)\cap {\rm supp}(g')=\emptyset$;
 \item [\it (c)]
 $\varphi(g) = \prod\limits_{p\in \mathbb{N}\diagup s}
 \varphi\left(s_p \gamma(p)\right)$ for each $g=s \gamma
 =\prod\limits_{p\in \mathbb{N}\diagup s}
 s_p \gamma(p)$ (see {\ref{decompositiontocycles}}).
\end{itemize}
 \end{Prop}
 \begin{proof} The equivalence of $(b)$ and $(c)$ is obvious. We prove the equivalence of $(a)$ and $(b)$.
   Using GNS-construction, we
 build the representation $\pi_{\varphi}$ of the group
$\Gamma\wr\mathfrak{S}_\infty$ which acts in the Hilbert space
$\mathcal{H}_\varphi$ with cyclic vector $\xi_\varphi$ such that
 \begin{eqnarray*}
 \varphi(g)=\left( \pi_\varphi\left( g\right)\xi_\varphi,\xi_\varphi
 \right)\textit{ for each }g\in\Gamma\wr\mathfrak{S}_\infty.
\end{eqnarray*}
 Suppose that the property ({\it a}) holds.
  Consider two elements
 $g=s  \gamma$ and
 $g^{\prime}=s^{\prime}  \gamma^{\prime}$
from $ \Gamma\wr\mathfrak{S}_\infty$ satisfying ${\rm supp}(g)\cap
{\rm supp}(g')=\emptyset$.
 Then there exists a sequence
 $\left\{ s_n \right\}_{n\in\mathbb{N}}\subset
 \mathfrak{S}_\infty$ such that for each $n$
 \begin{eqnarray}\label{assympt}
{\rm supp}(s_n)\cap {\rm supp}(g)=\emptyset\text{ and }
 {\rm supp}(s_n  g^{\prime}s_n^{-1})\subset\{n+1,n+2,\ldots\}.
  \end{eqnarray} For example we can put $s_n=\prod\limits_{i\in
  {\rm supp}(g')}(i,i+k+n)$, where $k$ is fixed number such that
  ${\rm supp}(g)\cup {\rm supp}(g')\subset\{1,2,\ldots,k\}$.
Using the ideas of the proof of the lemma \ref{weak-lim} we obtain,
that the limit $\lim\limits_{n\rightarrow
\infty}\pi_\varphi(s_ng's_n)$ exists in the weak operator topology
and the next equality holds:
 \begin{eqnarray}\label{fi s_n}
w-\lim\limits_{n\rightarrow \infty}\pi_\varphi(s_ng's_n)=\varphi(g')
I.
  \end{eqnarray}
Using ({\ref{assympt}}), (\ref{fi s_n}) and
$\mathfrak{S}_\infty$-centrality of $\varphi$,
 we obtain
 \begin{eqnarray*}
 &\varphi\left( g  g^\prime \right)=
 \lim\limits_{n\to\infty}\varphi
 \left( g  s_ng^\prime  s_n^{-1}\right)=\\&
  \lim\limits_{n\to\infty}\left(\pi_\varphi
 (g)\pi_\varphi\left(s_n g^\prime s_n^{-1}\right)\xi_\varphi,\xi_\varphi\right)= \varphi (g)  \varphi
\left(g^\prime\right).
 \end{eqnarray*}
 Thus {\it (b)} follows from {\it(a)}.

 Further suppose that the condition {\it(b)} holds.
      If $\pi_{\varphi}
\left(\Gamma\wr\mathfrak{S}_\infty\right)^\prime
\bigcap\pi_{\varphi}
\left(\Gamma\wr\mathfrak{S}_\infty\right)^{\prime\prime}=
\mathcal{Z}$
 is larger than the scalars, then it contains a pair
of
 orthogonal projections $E$ and $F$ satisfying the condition:
 \begin{eqnarray}\label{0inequalities}
  E  F=0.
 \end{eqnarray}
 Fix arbitrary $\varepsilon>0$. By the von Neumann Double Commutant Theorem there
exist $g_k, h_k\in\Gamma\wr\mathfrak{S}_\infty$ and complex numbers
 $c_k,d_k$ $\left( k=1,2,\ldots, N<\infty \right)$
such that
\begin{eqnarray}\label{inequalities}
\begin{split}
\left|\left|\sum\limits_{k=1}^Nc_k\pi_\varphi\left(g_k\right)\xi_\varphi
-E\xi_\varphi\right|\right|<\varepsilon,
\\ \left|\left|\sum\limits_{k=1}^N d_k\pi_\varphi\left(h_k\right)
\xi_\varphi-F\xi_\varphi\right|\right|<\varepsilon.
\end{split}
\end{eqnarray}
 Fix $n$ such that ${\rm supp}(g_k)\subset\{1,2,\ldots,n\}$ and
 ${\rm supp}(h_k)\subset\{1,2,\ldots,n\}$ for each $k$.
As $\varphi$ is $\mathfrak{S}_\infty$-central, using
({\ref{inequalities}}), we obtain
\begin{eqnarray}\label{1inequalities}
\left|\left|\sum\limits_{k=1}^N c_k\pi_\varphi\left(\omega_n
g_k\omega_n\right)\xi_\varphi-E\xi_\varphi\right|\right|<\varepsilon,\textit{
(see (\ref{omega}))}.
\end{eqnarray}
Now, using ({\ref{0inequalities}}), ({\ref{inequalities}}) and
 ({\ref{1inequalities}}), we have
\begin{eqnarray}\label{2inequalities}\left|\left(\sum\limits_{k=1}^Nc_k
\pi_\varphi\left(\omega_n g_k\omega_n \right) \sum\limits_{k=1}^N
d_k\pi_\varphi \left( h_k
\right)\xi_\varphi,\xi_\varphi\right)\right|<2\varepsilon+\varepsilon^2.
\end{eqnarray} Note, that ${\rm supp}(\omega_ng_k\omega_n)\subset\{n+1,n+2,\ldots\}$
for each $k$. Therefore, by the property $(b)$,
({\ref{inequalities}}) and
 ({\ref{1inequalities}}), one has:
\begin{eqnarray}\label{3inequalities}
\begin{split}\left|\left(\sum\limits_{k=1}^Nc_k\pi_\varphi\left(\omega_n g_k\omega_n
\right)  \sum\limits_{k=1}^N d_k\pi_\varphi \left( h_k
\right)\xi_\varphi,\xi_\varphi\right)\right|=\\\left|\left(\sum\limits_{k=1}^Nc_k\pi_\varphi\left(\omega_n
g_k\omega_n\right)\xi_\varphi,\xi_\varphi\right) \left(
\sum\limits_{k=1}^N d_k\pi_\varphi \left( h_k
\right)\xi_\varphi,\xi_\varphi\right)\right|>\\
\left(E\xi_\varphi,\xi_\varphi\right)
\left(F\xi_\varphi,\xi_\varphi\right)-\varepsilon\left(\left(E\xi_\varphi,\xi_\varphi\right)
+\left(F\xi_\varphi,\xi_\varphi\right)\right)-\varepsilon^2.\end{split}\end{eqnarray}
Note that, as $\xi_\varphi$ is cyclic, $E\xi_\varphi\neq 0$ and
$F\xi_\varphi\neq 0$. Therefore, taking in view
(\ref{2inequalities}) and (\ref{3inequalities}),
 we arrive at a
contradiction.
 \end{proof}
 Denote the element $ \sigma_n \in \mathfrak{S}_\infty$  by the formula:
 \begin{eqnarray}\label{12...n}
  \sigma_n (i)= \left\{
 \begin{array}{rl}
 i+1&\textit{ if }\; i<n,\\
 1&\textit{ if }\;i=n,\\
 i&\textit{ if }\;i>n.
 \end{array}\right.
\end{eqnarray}
 \begin{Co}\label{Co of mult}
 Each indecomposable $\mathfrak{S}_\infty$-central state
$\varphi$ on the group
  $ \Gamma\wr\mathfrak{S}_\infty$ is defined by its values on the
  elements of the form
  $ \sigma_n \gamma,$ where $\gamma=(\gamma_1,\gamma_2,\ldots,\gamma_n,e,e,\ldots)$ and $n\in
  \mathbb{N}$.
 \end{Co}
 \begin{proof}
  By the proposition \ref{multiplicativity}, $\varphi$ is defined by
  its values on the elements of the view $s_p\gamma(p)$ (see
  (\ref{decompositiontocycles})). Fix an element $s_p\gamma(p)$. Let
   $n=|p|$.
  Then there exists a permutation $h\in
  \mathfrak{S}_\infty$ such that $hs_ph^{-1}= \sigma_n $.
  Therefore
  $\varphi(s_p\gamma(p))=\varphi(hs_p\gamma(p)h^{-1})=\varphi( \sigma_n h\gamma(p)h^{-1})$,
  which proves the corollary.
 \end{proof}
 \paragraph{The characters of the group $\mathfrak{S}_\infty$ and $\Gamma\wr
\mathfrak{S}_\infty$.}
 In the paper {\cite{Thoma}}, E.Thoma obtained the following
 remarkable description of all {\it indecomposable}
 character ($\mathfrak{S}_\infty$-central states) of the group $\mathfrak{S}_\infty$. Characters of the group $\mathfrak{S}_\infty$ are labeled by a
 pair of non-increasing positive sequences of
 numbers
 $\left\{ \alpha_k \right\}$,
$\left\{ \beta_k \right\}$ $\left(k\in\mathbb{N} \right)$,
 such that
\begin{eqnarray}\label{cond}
\sum\limits_{k=1}^{\infty}\alpha_k +
 \sum\limits_{k=1}^{\infty}\beta_k\leq 1.
\end{eqnarray}
 The value of the corresponding character on
 a cycle of length $l$ is
\begin{eqnarray*}\label{color2}
\sum\limits_{k=1}^{\infty}\alpha_k^l +(-1)^{l-1}
 \sum\limits_{k=1}^{\infty}\beta_k^l.
\end{eqnarray*}
Its value on a product of several
 disjoint cycles equals to the product
of values on each of cycles.

 In $\cite{DN}$ authors described all indecomposable characters on
 the group $\Gamma\wr
\mathfrak{S}_\infty.$ Before to formulate the main result of
$\cite{DN}$ we introduce some more notations. We call an element
$g=s  \gamma$ a generated cycle if either $s$ is a cycle and
 ${\rm supp}(\gamma)\subset {\rm supp}(s)$ or $s=e$ and ${\rm supp}(\gamma)=\{n\}$ for some
 $n$.
For an element $g=s  \gamma$ and an orbit $p\in
 \mathbb{N}/s$ choose the minimal number $k\in p$ and denote
 \begin{eqnarray}\label{tilde gamma}\tilde{\gamma}(p)= \gamma_k \gamma_{s^{(-1)}(k)} \cdots
\gamma_{s^{(-l)}(k)}\cdots
   \gamma_{s^{(-|p|+1)}(k)}.
   \end{eqnarray}

  For a factor-representation $\tau$ of the finite type let $\chi_{\tau}$
  be its normalized character. That is  $\chi_\tau(g) =tr_{\mathcal{M}_\tau}\left( \tau(g)
\right)$, where $ tr_\mathcal{M}$  stands for the
 unique normal, normalized $\left( tr_\mathcal{M}(I)=1 \right)$
trace  on the factor $\mathcal{M}$ of the finite type. Note that
 $\chi_\tau(e)=1$. Let $tr$ be the ordinary matrix  normalized trace.
 \begin{Th}[\cite{DN}, \cite{DN1}]\label{mainth}
  Let $\varphi$ be a function on the group $\Gamma\wr
\mathfrak{S}_\infty$. Then the following conditions are equivalent.

$a)$ $\varphi$ is an indecomposable character.

$b)$ There exist a representation
 $\tau$ of the {\it finite} type of the group $\Gamma$,
 two  non-increasing positive sequences of  numbers
 $\left\{ \alpha_k \right\}$,
$\left\{ \beta_k \right\}$ $\left(k\in\mathbb{N} \right)$ and  two
sequences $\left\{\rho_k\right\},\left\{\varrho_k\right\}$ of
 finite-dimensional irreducible   representations of
 $\Gamma$ with properties
\begin{itemize}
\item {\rm(i)}
 $\delta=1-\sum\limits_k\alpha_k
dim\,\rho_k-\sum\limits_k\beta_k dim \varrho_k\geqslant 0$;
 \item {\rm(ii)} if $s$ is cycle, $g=s\gamma$ $\left(
\gamma \in\Gamma^\infty_e\right)$, $p={\rm supp} s ={\rm supp}
\left(s\gamma \right)$, then
\begin{eqnarray*}\varphi(g)=\left\{
 \begin{array}{ll}
\sum\limits_k\alpha_k\,tr(\rho_k(\gamma_n))+
\sum\limits_k\beta_k,\,tr(\varrho_k(\gamma_n))+\delta\chi_\tau(\gamma_n),
\textit{ if }\;p=\{n\},\\
\sum\limits_k\alpha_k^{|p|} tr(\rho_k(\tilde{\gamma}(p)))+(-1)^{|p\,|-1}
\sum\limits_k\beta_k^{|p\,|} tr(\varrho_k(\tilde{\gamma}(p))),
 \textit{
if }\;|p\,|>1;\end{array}\right.
\end{eqnarray*}
\item {\rm(iii)} if $g=s \gamma
 =\prod\limits_{p\in \mathbb{N}\diagup s}
 s_p \gamma(p)$ (see {\ref{decompositiontocycles}}), then
 $\varphi(g) = \prod\limits_{p\in \mathbb{N}\diagup s}
 \varphi\left(s_p \gamma(p)\right)$.
\end{itemize}
 \end{Th}


\subsubsection{Examples of representations.}\label{ex of repr}
\paragraph{ Parameters of states.}\label{paragraph2.1}
Let $A$ be  a self-adjoint operator of the {\it trace class} (see
{\cite{RS}}) from $\mathcal{B}(\mathcal{H})$ with the property:

 ${\rm Tr} (|A|)\leq1$, where ${\rm Tr}$ is ordinary trace\footnote{If
      $\mathfrak{p}$ is the minimal projection from
      $\mathcal{B}(\mathcal{H})$, then ${\rm Tr} (\mathfrak{p})=1$.} on
      $\mathcal{B}(\mathcal{H})$.

 Further we fix vector $\hat{\xi}\in{\rm Ker}
     \, A$ and the unitary representation $\rho$ of $\Gamma$ in $\mathcal{H}$,
which satisfies the conditions:
\begin{itemize}
  \item {\rm(1)} if ${\rm Tr} (|A|)=1$, then subspace $({\rm Ker}
   \,   A)^\perp=\mathcal{H}\ominus{\rm Ker}\, A$  is cyclic for
      $w^*$-algebra $\mathfrak{A}$ generated by $A$ and
      $\rho(\Gamma)$;
  \item {\rm(2)} if ${\rm Tr} (|A|)<1$,
   subspace $\widetilde{\mathcal{H}}$ is
      generated by $\left\{ \mathfrak{A}v,   v\in ({\rm
      Ker} A)^\perp \right\}$  and $\mathcal{H}_{reg}=
      \mathcal{H}\ominus\widetilde{\mathcal{H}}$, then
 $\dim \mathcal{H}_{reg}=\infty$;
  \item {\rm(3)} \label{page11} if $P_{]0,1]}$ and $P_{[-1,0[}$ are the spectral
  projections of $A$, then subspaces $\mathcal{H}_+$ and
  $\mathcal{H}_-$ generated by vectors
  $ \left\{\mathfrak{A}v,\;v\in P_{]0,1]}\mathcal{H} \right\}$ and
  $ \left\{\mathfrak{A}v,\;v\in P_{[-1,0[}\mathcal{H}
  \right\}$, respectively, are orthogonal;
  \item {\rm(4)} \label{property4}there exist ${\rm I}_\infty$-factor
      $N^\prime_{reg}\subset \left(\rho\left(\Gamma
      \right)\Big|_{\mathcal{H}_{reg}}\right)^\prime$ with matrix unit
      $ \left\{\mathfrak{e}_{kl}^\prime,\; k,l\in\mathbb{N} \right\}$
      such that  $\hat{\xi}\in
      \mathfrak{e}_{11}^\prime\mathcal{H}_{reg}$,
      $\left\|\hat{\xi}\right\|=1$ and
      $\mathfrak{e}_{11}^\prime\mathcal{H}_{reg}$ is generated by $
      \left\{\rho \left(\Gamma \right)\hat{\xi} \right\}$. In
      particular, if ${\rm Tr} (|A|)=1$ then $\hat{\xi}=0$.
      When  ${\rm Tr} (|A|)<1$ we assume for convenience that
      $\left\|\hat{\xi}\right\|=1$.
\end{itemize}
\paragraph{Hilbert space $\mathcal{H}^\rho_A$.}\label{paragraph2.2}
 Define
a state $\psi_k$ on $\mathcal{B}\left(\mathcal{H} \right)$ as
follows
\begin{eqnarray}\label{psik}
\psi_k\left(v \right)={\rm Tr}\left(v|A| \right)+\left(1-{\rm
Tr}\left(|A| \right)
\right)\left(v\mathfrak{e}_{k1}^\prime\hat{\xi},\mathfrak{e}_{k1}^\prime\hat{\xi}
\right),\;\;v\in\mathcal{B}\left(\mathcal{H} \right).
\end{eqnarray}
Let $ _1\psi_k$ denote the product-state on $\mathcal{B} \left(
H\right)^{\otimes k}$:
\begin{eqnarray}\label{psik1}
 _1\psi_k\left(v_1\otimes v_2\otimes \ldots\otimes v_k
 \right)=\prod\limits_{j=1}^k\psi_j\left(v_j \right).
\end{eqnarray}
Now define inner product on $\mathcal{B} \left( H\right)^{\otimes
k}$ by
\begin{eqnarray}\label{psik2}
\left( v,u\right)_k=\,_1\psi_k\left( u^*v \right).
\end{eqnarray}
Let $\mathcal{H}_k$ denote the Hilbert space obtained by completing
$\mathcal{B} \left( H\right)^{\otimes k}$  in above inner product
norm. Now we consider the natural isometrical embedding
\begin{eqnarray}
v\ni\mathcal{H}_k\mapsto v\otimes {\rm I}\in\mathcal{H}_{k+1}.
\end{eqnarray}
and define Hilbert space $\mathcal{H}^\rho_A$ as completing
$\bigcup\limits_{k=1}^\infty \mathcal{H}_k$.
\paragraph{The action $\Gamma\wr\mathfrak{S}_{\infty}$ on
$\mathcal{H}^\rho_A$.}\label{paragraph2.3}
  First, using the
embedding $a\ni\mathcal{B}\left(\mathcal{H} \right)^{\otimes
k}\mapsto a\otimes{\rm I}\in\mathcal{B}\left(\mathcal{H}
\right)^{\otimes k+1}$, we identify $\mathcal{B}\left(\mathcal{H}
\right)^{\otimes k}$ with subalgebra $\mathcal{B}\left(\mathcal{H}
\right)^{\otimes
k}\otimes\mathbb{C}\subset\mathcal{B}\left(\mathcal{H}
\right)^{\otimes k+1}$. Therefore, algebra
$\mathcal{B}\left(\mathcal{H}
\right)^{\otimes\infty}=\bigcup\limits_{n=1}^\infty
\mathcal{B}\left(\mathcal{H} \right)^{\otimes n}$ is well defined.

Further we give the explicit embedding $\mathfrak{S}_\infty$ into
unitary group of $\mathcal{B}\left(\mathcal{H}
\right)^{\otimes\infty}$. First fix the matrix unit $ \left\{e_{pq}:
p,q=1,2,\ldots,n={\rm
dim}\,\mathcal{H}\right\}\subset\mathcal{B}\left(\mathcal{H}
\right)$ with the properties:
\begin{itemize}
    \item {\rm (i)} projection $e_{kk}$ is minimal and
    $e_{kk}A=c_{kk}e_{kk}$ $\left(c_{kk}\in\mathbb{C} \right)$ for all
    $k=1,2,\ldots,n$;
    \item {\rm (ii)} $e_{kk}\mathcal{H}_+\subset\mathcal{H}_+$ and
    $e_{kk}\mathcal{H}_-\subset\mathcal{H}_-$  for all
    $k=1,2,\ldots,n$.
\end{itemize}
Put $X=  \left\{1,2,\ldots,n \right\}^{\times\infty}$. For
$x=\left(x_1,x_2,\ldots,x_l,\ldots \right)\in X$  we set
$\mathfrak{l}_A(x)=\left|
\left\{i:e_{x_i\,x_i}\mathcal{H}\subset\mathcal{H}_- \right\}
\right|$. Define subsequence $x_A=\left(x_{i_1},x_{i_2},\ldots
x_{i_l},\ldots \right)\in
 \left\{ 1,2,\ldots,n\right\}^{\mathfrak{l}_A(x)}$ by induction
\begin{eqnarray}
i_1= {\rm min}\left\{i:e_{x_i\,x_i}\mathcal{H}\subset\mathcal{H}_-
\right\} \text{ and }\;i_k={\rm min}
\left\{i>i_{k-1}:e_{x_i\,x_i}\mathcal{H}\subset\mathcal{H}_-
\right\}.
\end{eqnarray}
 For $s\in\mathfrak{S}_\infty$ denote by $c(x,s)$ the unique
 permutation from
$\mathfrak{S}_{{\mathfrak{l}_A(x)}}\subset\mathfrak{S}_\infty$ such that
 \begin{eqnarray}
s^{-1}\left(i_{c(x,s)(1)} \right)<s^{-1}\left(i_{c(x,s)(2)}
\right)<\ldots<s^{-1}\left(i_{c(x,s)(l)} \right)<\ldots ..
 \end{eqnarray}
Let $\mathfrak{S}_\infty$ acts on $X$ as follows
\begin{eqnarray}
X\times\mathfrak{S}_\infty\ni(x,s)\mapsto
sx=\left(x_{s(1)},x_{s(2)},\ldots,x_{s(l)},\ldots \right)\in X.
\end{eqnarray}
By definition, $\left(
sx\right)_A=\left(x_{i_{c(x,s)(1)}},x_{i_{c(x,s)(2)}},\ldots,x_{i_{c(x,s)(l)}},\ldots
\right)$. Therefore,
\begin{eqnarray}\label{cocycle}
c(x,ts)=c(sx,t)c(x,s)\;\; \text{ for all
}\;t,s\in\mathfrak{S}_\infty;x\in X.
\end{eqnarray}
Given any $s\in\mathfrak{S}_\infty$ put
\begin{eqnarray*}
U_N(s)=\sum\limits_{ x_1,x_2, \ldots,x_N=1}^n {\rm sign}\,
\left(c(x,s)\right)\;e_{x_{s(1)}\,x_1}\otimes
e_{x_{s(2)}\,x_2}\otimes\ldots\otimes e_{x_{s(N)}\,x_N},
\end{eqnarray*}
where $N<\infty$ satisfies the condition: $s(i)=i$ for all $i\geq
N$, $x=\left(x_1,x_2,\ldots,x_N,\ldots \right)$.
 We see at once that for $L>N$
 \begin{eqnarray*}
U_N(s)\otimes\underbrace{{\rm I}\otimes{\rm
I}\otimes\ldots\otimes{\rm I}}_{L-N}=U_L(s).
 \end{eqnarray*}
Thus operator $U(s)=U_N(s)\otimes{\rm I}\otimes{\rm
I}\otimes\ldots\in \mathcal{B}\left(\mathcal{H}
\right)^{\otimes\infty}=\bigcup\limits_{n=1}^\infty
\mathcal{B}\left(\mathcal{H} \right)^{\otimes n}$ is well defined.
 It  follows from
\ref{cocycle} that
\begin{eqnarray}
U(t)U(s)=U(ts) \text{ for all } t,s\in\mathfrak{S}_\infty.
\end{eqnarray}
It is clear that
\begin{eqnarray*}
{\rm sign}\,\left(c(x,s)c(y,s) \right)U(s)\left(e_{x_1\,y_1}\otimes
e_{x_2\,y_2}\otimes\ldots\otimes e_{x_N\,y_N}\otimes {\rm
I}\otimes{\rm I}\ldots\otimes{\rm
I}\otimes\ldots \right)U(s)^*\\
 =e_{x_{s^{-1}(1)}\,y_{s^{-1}(1)}}\otimes e_{x_{s^{-1}(2)}\,y_{s^{-1}(2)}}\otimes\ldots\otimes
e_{x_{s^{-1}(N)}\,y_{s^{-1}(N)}}\otimes {\rm I}\otimes{\rm
I}\ldots\otimes{\rm I}\otimes\ldots.
\end{eqnarray*}
If $x$, and $y$ satisfies the condition:

$e_{x_i\,x_i}\mathcal{H}\subset\mathcal{H}_-$ if and only if, when
$e_{y_i\,y_i}\mathcal{H}\subset\mathcal{H}_-$,

\noindent then, by definition cocycle $c$, we have $c(x,s)=c(y,s)$.
Therefore,
\begin{eqnarray}\label{UmatrixUnit}
\begin{split}
U(s)\left(e_{x_1\,y_1}\otimes e_{x_2\,y_2}\otimes\ldots\otimes
e_{x_N\,y_N}\otimes {\rm I}\otimes{\rm I}\ldots \right)U(s)^*\\
 =e_{x_{s^{-1}(1)}\,y_{s^{-1}(1)}}\otimes e_{x_{s^{-1}(2)}\,y_{s^{-1}(2)}}\otimes\ldots\otimes
e_{x_{s^{-1}(N)}\,y_{s^{-1}(N)}}\otimes {\rm I}\otimes{\rm I}\ldots.
\end{split}
\end{eqnarray}
Hence, using properties {\rm(2)}-{\rm(3)} on the page
\pageref{page11}, we obtain
\begin{eqnarray}\label{relurho}
\begin{split}
U(s)\left(\rho\left(\gamma_1\right)\otimes\rho\left(\gamma_2\right)\otimes\ldots
\otimes\rho\left(\gamma_N\right)\otimes \ldots
\right)U(s)^*\\
=\rho\left(\gamma_{s^{-1}(1)}\right)\otimes
 \rho\left(\gamma_{s^{-1}(2)}\right)\otimes
\ldots\otimes\rho\left(\gamma_{s^{-1}(N)}\right)\otimes\ldots
\end{split}
\end{eqnarray}
for all $s\in\mathfrak{S}_\infty$, $\gamma_l\in\Gamma $.

Now we define the operators  $\Pi_A^\rho(s)$,
$\left(s\in\mathfrak{S}_\infty \right)$ and $\Pi_A^\rho(\gamma )$,
$\left(\gamma =\left(\gamma_1,\gamma_2,\ldots
\right)\in\Gamma^\infty_e \right)$ on $\mathcal{H}_A^\rho$ as
follows
\begin{eqnarray}\label{Piarho}
\begin{split}
\Pi_A^\rho(s)v=U(s)v,\;\;\;v\in \mathcal{H}_A^\rho;\\
\Pi_A^\rho(\gamma)v=\left(\rho\left(\gamma_1\right)\otimes\rho\left(\gamma_2\right)\otimes\ldots
 \right)v.
\end{split}
\end{eqnarray}
By (\ref{relurho}), $\Pi_A^\rho$ can be extended to the unitary
representation of $\Gamma\wr\mathfrak{S}_\infty$.

The next proposition follows from the definition of Hilbert space
$\mathcal{H}_A^\rho$ (see paragraph \ref{paragraph2.2}) and
proposition \ref{multiplicativity}.
\begin{Prop}\label{Prop11a}
Let $I$ be the unit in $\mathcal{B}\left(\mathcal{H}
\right)^{\otimes\infty}$. Identify the elements of
$\mathcal{B}\left(\mathcal{H} \right)^{\otimes\infty}$ with the
corresponding vectors in $\mathcal{H}_A^\rho$. Put $\psi_A^\rho
\left(s\gamma \right)=\left(\Pi_A^\rho(s)\Pi_A^\rho(\gamma )I,I
\right)$. Then $\phi_A^\rho$ is indecomposable
$\mathfrak{S}_\infty$-central state on
$\Gamma\wr\mathfrak{S}_\infty$ (see definitions \ref{indecomposable}
and \ref{def central}).
\end{Prop}
Let $A_1$, $A_2$ be  the self-adjoint operators  of the {\it trace
class} (see {\cite{RS}}) from $\mathcal{B}(\mathcal{H})$ with the
property ${\rm Tr} (|A_j|)\leq1$, $\left(j=1,2 \right)$, and let
$\rho_1$, $\rho_2$ be the unitary representations of $\Gamma$:
$\rho_i:\gamma\in \Gamma\mapsto\rho_i(\gamma)\in
\mathcal{B}(\mathcal{H})$.
\begin{Prop}\label{Prop12}
Let $\left(\mathcal{H}_i,A_i,\rho_i,\hat{\xi}_i\right)$, $i=1,2$
satisfy assumptions {\rm(1)}-{\rm(4)} (paragraph
\ref{paragraph2.1}). Equality
$\psi_{A_1}^{\rho_1}=\psi_{A_2}^{\rho_2}$ holds if and only if there
exists isometry  $\mathcal{U}:\mathcal{H}_1\mapsto\mathcal{H}_2$ such that
\begin{eqnarray}\label{unitaryeq}
\hat{\xi}_2=\mathcal{U}\hat{\xi}_1,\;A_2=\mathcal{U}A_1\mathcal{U}^{-1}
\text{ and }\;\; \rho_2(\gamma )=\mathcal{U}\rho_1(\gamma
)\mathcal{U}^{-1} \text{ for all } \gamma \in\Gamma.
\end{eqnarray}
\end{Prop}
\begin{proof}
Assume (\ref{unitaryeq}) hold. It follows from (\ref{psik}) and
proposition \ref{Prop11a} that
$\psi_{A_1}^{\rho_1}=\psi_{A_2}^{\rho_2}$.

Conversely, suppose that $\psi_{A_1}^{\rho_1}=\psi_{A_2}^{\rho_2}$.

 Denote by $\Pi_A^{\rho\,0}$  the restriction $\Pi_A^\rho$ to subspace $
\left[ \Pi_A^\rho\left(\Gamma\wr\mathfrak{S}_\infty \right)I\right]$
generated by the vectors $
\left\{\Pi_A^\rho\left(\Gamma\wr\mathfrak{S}_\infty \right)I
\right\}$. Let $\left(l\,k \right)$ be the transposition
interchanging $l$ and $k$. According to the construction of
representation $\Pi_A^\rho$ and properties {\rm (i)-(ii)} from
paragraph \ref{paragraph2.3}, there exists operator
\begin{eqnarray}\label{astrans}
\mathcal{O}_l=w-\lim\limits_{k\to\infty}\Pi_A^\rho\left( \left(l\,k
\right)\right)
\end{eqnarray}
and
\begin{eqnarray}\label{Oaction}
\mathcal{O}_l\left(a_1\otimes a_2\otimes\ldots\right)=b_1\otimes
b_2\otimes\ldots,\text{ where }
 b_k=\left\{\begin{array}{ll}a_k,&\textit{ if }\;k\neq l,\\
Aa_k,&\textit{ if }\;k=l. \end{array}\right.
\end{eqnarray}
Let $\mathfrak{A}_l^{A\,\rho}$ be $w^*$-algebra in
$\Pi_A^\rho\left(\Gamma\wr\mathfrak{S}_\infty
\right)^{\prime\prime}$ generated by $\mathcal{O}_l$ and
$\underbrace{{\rm I}\otimes\ldots\otimes{\rm
I}}_{l-1}\otimes\rho(\gamma )\otimes{\rm I}\otimes{\rm
I}\otimes\ldots$, $\gamma \in\Gamma$. Denote by $\mathcal{P}_0$ the
orthogonal projection $\mathcal{H}_A^\rho$ onto $ \left[
\mathfrak{A}_l^{A\,\rho} I\right]$.

 First we prove that $w^*$-algebra $
\left\{A,\rho(\Gamma)
\right\}^{\prime\prime}\subset\mathcal{B}(\mathcal{H})$ generated by
$A$ and $\rho(\Gamma)$ is isomorphic to $w^*$-algebra
$\mathfrak{A}_l^{A\,\rho}\mathcal{P}_0$. Namely, the map
\begin{eqnarray}\begin{split}
\mathfrak{m}_l:A\mapsto \mathcal{O}_l\mathcal{P}_0,\\
\mathfrak{m}_l:\rho(\gamma)\mapsto \left(\underbrace{{\rm
I}\otimes\ldots\otimes{\rm I}}_{l-1}\otimes\rho(\gamma )\otimes{\rm
I}\otimes{\rm I}\otimes\ldots\right)\mathcal{P}_0
\end{split}\end{eqnarray} extends to an isomorphism of $\left\{A,\rho(\Gamma)
\right\}^{\prime\prime}$ onto
$\mathfrak{A}_l^{A\,\rho}\mathcal{P}_0$.

Using (\ref{Oaction}) and definition of $\Pi_A^\rho$, we can
consider $\mathfrak{m}_l$ as the GNS-representation of $
\left\{A,\rho(\Gamma)
\right\}^{\prime\prime}\subset\mathcal{B}(\mathcal{H})$
corresponding to $\psi_k$ (see (\ref{psik})).  Thus ${\rm Ker}\,
\mathfrak{m}_l= \left\{a\in \left\{A,\rho(\Gamma)
\right\}^{\prime\prime}:\mathfrak{m}_l\,(a)=0\right\}$ is weakly
closed two-sided ideal. Therefore, there exists unique orthogonal
projection $e$ from the center of $\left\{A,\rho(\Gamma)
\right\}^{\prime\prime}$ such that
\begin{eqnarray}\label{Kerml}
{\rm Ker}\,\mathfrak{m}_l
=e\left\{A,\rho(\Gamma)\right\}^{\prime\prime} \text{ (see
\cite{Tak1})}.
\end{eqnarray}
Let us prove that $e=0$.

Denote by $c\left(\widetilde{P}\right)$ central support of
orthogonal projection $\widetilde{P}\in \left\{A,\rho(\Gamma)
\right\}^\prime$: $\widetilde{P}\mathcal{H}=\widetilde{\mathcal{H}}$
(see property {\rm (2)} from paragraph \ref{paragraph2.1}).

Let us first show that
\begin{eqnarray}\label{ecentrsupp}
e\,c\left(\widetilde{P}\right)=0.
\end{eqnarray}
 Conversely, suppose that $e\,c\left(\widetilde{P}
\right)\neq 0$. Hence, since the map $\left\{A,\rho(\Gamma)
\right\}^{\prime\prime}\,c\left(\widetilde{P}\right)\ni a\mapsto
a\widetilde{P}\in\left\{A,\rho(\Gamma)
\right\}^{\prime\prime}\widetilde{P}$ is isomorphism, we obtain
$e\,\widetilde{P}\neq0$. It follows from properties {\rm (1)}-{\rm
(3)} (paragraph \ref{paragraph2.1})) that
 $e\,\left(P_{]0,1]}+P_{[-1,0[}\right)\neq 0$. Thus, by
 (\ref{psik}), $\psi_l(e)\neq 0$. Therefore,
 $e\notin {\rm Ker}\,\mathfrak{m}_l$. This
 contradicts property (\ref{Kerml}).

Now, using (\ref{ecentrsupp}) and property {\rm (2)} (paragraph
\ref{paragraph2.1})), we have
\begin{eqnarray}
e\,\left(I-c\left(\widetilde{P} \right) \right)\mathcal{H}\subseteq
\mathcal{H}_{reg}.
\end{eqnarray}
Therefore, if $e\,\left(I-c\left(\widetilde{P} \right)
\right)\neq0$, then, using property {\rm (4)} (paragraph
\ref{paragraph2.1}), we obtain
\begin{eqnarray}
e\,\left(I-c\left(\widetilde{P} \right)
\right)\mathfrak{e}_{l1}^\prime\hat{\xi }\neq0.
\end{eqnarray}
Again,  by
 (\ref{psik}), $\psi_l(e)\neq 0$ and $e\notin {\rm
 Ker}\,\mathfrak{m}_l$. It follows from (\ref{Kerml}) that
 \begin{eqnarray}
e\,\left(I-c\left(\widetilde{P} \right) \right)=0.
\end{eqnarray}
Hence, using (\ref{ecentrsupp}), we obtain
\begin{eqnarray}\label{Kerzero}
{\rm Ker}\,\mathfrak{m}_l =0.
\end{eqnarray}
Now we suppose that $\phi_{A_1}^{\rho_1}=\phi_{A_2}^{\rho_2}$. Let
$\mathcal{O}_l^{(1)}$ and $\mathcal{O}_l^{(2)}$ be the operators,
which are defined by  formula (\ref{astrans})  for representations
$\Pi_{A_1}^{\rho_1}$ and  $\Pi_{A_2}^{\rho_2}$ respectively. If
$\mathfrak{I}_l$ is the extension the map
\begin{eqnarray*}
\mathcal{O}_l^{(1)}\mathcal{P}_0\mapsto\mathcal{O}_l^{(2)}\mathcal{P}_0,\\
\underbrace{{\rm I}\otimes\ldots\otimes{\rm
I}}_{l-1}\otimes\rho_1(\gamma )\otimes{\rm I}\otimes{\rm
I}\otimes\ldots\mapsto \underbrace{{\rm I}\otimes\ldots\otimes{\rm
I}}_{l-1}\otimes\rho_2(\gamma )\otimes{\rm I}\otimes{\rm
I}\otimes\ldots.
\end{eqnarray*}
 by multiplication, then
 \begin{eqnarray}\label{invform}
\left(\mathfrak{I}_l(a)I,\mathfrak{I}_l(b)I\right)=\left(aI,bI\right)
\text{ for all } \;a,b\in \mathfrak{A}_l^{A_1\,\rho_1}\mathcal{P}_0.
 \end{eqnarray}
 It follows from (\ref{Kerzero}) that the map
 \begin{eqnarray}\label{thetaiso}
\left\{A_1,\rho_1(\Gamma) \right\}^{\prime\prime}\ni
a\stackrel{\theta}{\mapsto} \mathfrak{m}_l^{-1}
\circ\mathfrak{I}_l\circ\mathfrak{m}_l(a)\in
\left\{A_2,\rho_2(\Gamma) \right\}^{\prime\prime}
 \end{eqnarray}
 is an isomorphism. Since $\phi_{A_1}^{\rho_1}=\phi_{A_2}^{\rho_2}$,
 then, using definition of $\phi_{A}^{\rho}$, in particular (\ref{psik}), obtain for all
 $ v\in\left\{A_1,\rho_1(\Gamma) \right\}^{\prime\prime}$:
 \begin{eqnarray}\label{psieq}
 \begin{split}
{\rm Tr}\left(v|A_1| \right)+\left(1-{\rm Tr}\left(|A_1| \right)
\right)\left(v\hat{\xi}_1,\hat{\xi}_1 \right)\\= {\rm
Tr}\left(\theta(v)|A_2| \right)+\left(1-{\rm Tr}\left(|A_2| \right)
\right)\left(\theta(v)\hat{\xi}_2,\hat{\xi}_2 \right).
\end{split}
 \end{eqnarray}

 Without loss of generality we can assume that
 $\left\{A_1,\rho_1(\Gamma) \right\}^{\prime\prime},
 \left\{A_2,\rho_2(\Gamma) \right\}^{\prime\prime}\subset\mathcal{B}\left(\mathcal{H}
 \right)$. Let $P_{[-1,0[}^{(i)}$, $P_{]0,1]}^{(i)}$ be the spectral
 projections of $A_i$ $(i=1,2)$. Put
 $P_\pm^{(i)}=P_{[-1,0[}^{(i)}+P_{]0,1]}^{(i)}$. It is clear $\left({\rm Ker}\, A_i\right)^\perp=
 P_\pm^{(i)}\mathcal{H}$. Denote by $\widetilde{\mathcal{H}}_i$
  subspace $ \left[\left\{A_i,\rho_i(\Gamma) \right\}^{\prime\prime}
  P_\pm^{(i)}\mathcal{H}_i\right]$. Let $\widetilde{P}_i$ be the
  orthogonal projection of $\mathcal{H}_i$ onto
  $\widetilde{\mathcal{H}}_i$. Put $P^{(i)}_{reg}=I-\widetilde{P}_i$.  For $\alpha \in {\rm Spectrum }\, A_i$
  denote by $P_\alpha^{(i)}$ the corresponding spectral projection.

  Now, using properties of $\left(A_i,\rho_i \right)$ (see paragraph
  {\ref{paragraph2.1}}), we have
  \begin{eqnarray}
 {\rm
dim}\,P_\alpha^{(i)}\mathcal{H}<\infty\;\text{ and
}\;P_\pm^{(i)}=\sum\limits_{\alpha \in{\rm Spectrum }\, A_i:\alpha
\neq0}P_\alpha^{(i)}.
  \end{eqnarray}
  Therefore, there exists collection $ \left\{c^{(i)}_j
  \right\}_{j=1}^{N}$ of pairwise orthogonal projections from the
  center of $w^*$-algebra
   $P_\pm^{(i)}\left\{A_i,\rho_i(\Gamma)
   \right\}^{\prime\prime}P_\pm^{(i)}$ with properties
   \begin{eqnarray}\label{properties of central proj}
   \begin{split}
   \theta\left( c^{(1)}_j\right)=c^{(2)}_j \text{ (see (\ref{thetaiso})) };\;\;\;\sum\limits_{j=1}^{N}c^{(i)}_j=P_\pm^{(i)};\\
 c^{(i)}_jP_\pm^{(i)}\left\{A_i,\rho_i(\Gamma)
   \right\}^{\prime\prime}P_\pm^{(i)}c^{(i)}_j \text{ is a factor of type
   } I_{n_j}.
   \end{split}
   \end{eqnarray}
Fix matrix unit $ \left\{f^{(j)}_{k\,l} \right\}_{k,l=1}^{n_j}\subset
c^{(1)}_jP_\pm^{(1)}\left\{A_1,\rho_1(\Gamma)
   \right\}^{\prime\prime}P_\pm^{(1)}c^{(1)}_j$, which is a linear basis in
$c^{(1)}_jP_\pm^{(1)}\left\{A_1,\rho_1(\Gamma)
   \right\}^{\prime\prime}P_\pm^{(1)}c^{(1)}_j$, minimal projections
$\left\{f^{(j)}_{k\,k} \right\}_{k=1}^{n_j}$ satisfy condition
\begin{eqnarray}\label{spectrprmatrixunit}
P_\alpha^{(1)}f^{(j)}_{k\,k} =f^{(j)}_{k\,k}P_\alpha^{(1)} \;\text{
for all } \alpha \in {\rm Spectrum}\, A_1;\; k,j\in\mathbb{N}.
\end{eqnarray}
Now, using (\ref{invform}), (\ref{thetaiso}), (\ref{psieq}) and
definition of $\Pi_A^\rho$ (see paragraphs
\ref{paragraph2.1},\ref{paragraph2.2}, \ref{paragraph2.3}), we have
\begin{eqnarray}
{\rm Tr}\, \left( f^{(j)}_{k\,k}\right)= {\rm Tr}\, \left(\theta\left(
f^{(j)}_{k\,k}\right)\right) \text{ for all } \; k,j\in\mathbb{N}.
\end{eqnarray}
Therefore, there exists isometry $U:P_\pm^{(1)}\mathcal{H}_1\mapsto
P_\pm^{(1)}\mathcal{H}_2 $ such that
$UP_\pm^{(1)}\mathcal{H}_1=P_\pm^{(1)}\mathcal{H}_2$ and
\begin{eqnarray}\label{isomorphismSp}
Uf^{(j)}_{k\,k}U^{-1}=\theta\left( f^{(j)}_{k\,k}\right)
 \text{ for  } \; k=1,2,\ldots n_j;\;\;j=1,2,\ldots,N.
 \end{eqnarray}
Let $\mathcal{C}_i$ be the center of $w^*$-algebra
$\left\{A_i,\rho_i(\Gamma)
   \right\}^{\prime\prime}$ and let
   $c\left(P_\pm^{(i)}\right)\in\mathcal{C}_i$ be the central support
    of $P_\pm^{(i)}$. It follows from this and (\ref{properties of central proj})
that there exist pairwise orthogonal projections
 $ \left\{C^{(i)}_j \right\}_{j=1}^{N}\subset c\left(P_\pm^{(i)}\right)\cdot\mathcal{C}_i $
 with the next properties
 \begin{eqnarray}
 \begin{split}
c^{(i)}_j= C^{(i)}_j\cdot
P_\pm^{(i)},\;\;\;\;\sum\limits_{j=1}^NC^{(i)}_j=c\left(
P_\pm^{(i)}\right),\\
C^{(i)}_j\left\{A_i,\rho_i(\Gamma)
   \right\}^{\prime\prime}C^{(i)}_j \text{ is a factor of type
   } I_{N_j}.
   \end{split}
 \end{eqnarray}
In $C^{(1)}_j\left\{A_1,\rho_1(\Gamma)
   \right\}^{\prime\prime}C^{(1)}_j$ there exists matrix unit
   $ \left\{f^{(j)}_{k\,l} \right\}_{k,l=1}^{N_j}$ $\left(n_j\geq N_j \right)$.
Now, applying (\ref{isomorphismSp}), we obtain that
\begin{eqnarray}
\widetilde{U}=\sum_{j=1}^N\sum_{k=1}^{N_j}\theta\left(f_{k1}^{(j)}\right)Uf_{1k}
\end{eqnarray}
is an isometry of $c\left(P_\pm^{(1)}\right)\mathcal{H}_1 $ onto
$c\left(P_\pm^{(2)}\right)\mathcal{H}_2 $. An easy computation shows
that $\widetilde{U}f_{kl}^{(j)}
\widetilde{U}^{-1}=\theta\left(f_{kl}^{(j)} \right)$ for
$k,l=1,2,\ldots,N_j;$ $j=1,2,\ldots,N$. Thus
\begin{eqnarray}\label{Utilde}
\theta(a)=\widetilde{U}a \widetilde{U}^{-1}\; \text{ for all }\;
a\in c\left(P_\pm^{(1)}\right)\left\{A_1,\rho_1(\Gamma)
   \right\}^{\prime\prime}.
\end{eqnarray}
Hence, using (\ref{psieq}) and relations
$\theta\left(\left|A_1\right| \right)=\left|A_2\right|$,
$\theta\left(c\left(P_\pm^{(1)}\right)
\right)=c\left(P_\pm^{(2)}\right)$, which follows from the
definition of $\theta$ (see (\ref{thetaiso})), we have
\begin{eqnarray}\label{0.68}
\begin{split}
\left(\left(I-c\left(P_\pm^{(2)}\right)
\right)\theta(v)\hat{\xi}_2,\hat{\xi}_2
\right)=\left(\left(I-c\left(P_\pm^{(1)}\right)
\right)v\hat{\xi}_1,\hat{\xi}_1\right).
\end{split}
\end{eqnarray}
Since $\widetilde{P}_i\leq c\left(P_\pm^{(i)}\right)$, then
\begin{eqnarray}\label{0.69}
I-c\left(P_\pm^{(i)}\right)\leq P_{reg}^{(i)}, \;\;\; i=1,2.
\end{eqnarray}
Denote by
$\left\{\mathfrak{e}_{kl}^{(i)\prime},\;k,l\in\mathbb{N}\right\}$
  $\left(i=1,2 \right)$ the matrix unit from property {\rm (4)} of paragraph
\ref{paragraph2.1}. Now we define map $V$ as follows
\begin{eqnarray*}
 a\left(I-c\left(P_\pm^{(1)}\right)\right)\hat{\xi}_1
 \stackrel{V}{\mapsto}\theta(a) \left(I-c\left(P_\pm^{(2)}\right)\right)\hat{\xi}_2, \;\text{ where }\;
a\in\left\{A_1,\rho_1(\Gamma)
   \right\}^{\prime\prime}.
\end{eqnarray*}
By (\ref{0.68}) and  (\ref{0.68}), $V$ extends to  isometry $V$
of $\left(I-c\left(P_\pm^{(1)}\right)\right)\mathfrak{e}_{11}^{(1)\prime}
\mathcal{H}_1\subset P_{reg}^{(1)}\mathcal{H}_1$
onto $\left(I-c\left(P_\pm^{(2)}\right)\right)\mathfrak{e}_{11}^{(1)\prime}
\mathcal{H}_2\subset P_{reg}^{(2)}\mathcal{H}_2$ and for all
$a\in\left\{A_1,\rho_1(\Gamma)
   \right\}^{\prime\prime}$
\begin{eqnarray*}
V\left(I-c\left(P_\pm^{(1)}\right)\right)a\mathfrak{e}_{11}^{(1)\prime}V^{-1}=
\left(I-c\left(P_\pm^{(2)}\right)\right)\theta(a)\mathfrak{e}_{11}^{(2)\prime}.
\end{eqnarray*}
It follows from this that $\widetilde{V}=\sum\limits_{k=1}^\infty
 \mathfrak{e}_{k1}^{(2)\prime}V \left(I-c\left(P_\pm^{(1)}\right)\right)\mathfrak{e}_{1k}^{(1)\prime}$
is an isometry of $\left(I-c\left(P_\pm^{(1)}\right)\right)
\mathcal{H}_1$ onto $\left(I-c\left(P_\pm^{(2)}\right)\right)
\mathcal{H}_2$, satisfying the next relation
\begin{eqnarray*}
\widetilde{V}\left(I-c\left(P_\pm^{(1)}\right)\right)a\widetilde{V}^{-1}=
\left(I-c\left(P_\pm^{(2)}\right)\right)\theta(a) \;\;\;\;\; \left(a\in\left\{A_1,\rho_1(\Gamma)
   \right\}^{\prime\prime}\right).
\end{eqnarray*}
Hence, using (\ref{Utilde}), we obtain that
$W=\widetilde{U}c\left(P_\pm^{(1)}\right)+
\widetilde{V}\left(I-c\left(P_\pm^{(1)}\right)\right)$ is an isometry of
$\mathcal{H}_1$ onto $\mathcal{H}_2$
 and
\begin{eqnarray}
WaW^{-1}=\theta(a) \text{ for all }\;\; a\in\left\{A_1,\rho_1(\Gamma)
   \right\}^{\prime\prime}.
\end{eqnarray}
Now, on account of definition of $\theta$ and (\ref{psieq}) one can easy to check that
\begin{eqnarray}\label{0.71}
\begin{split}
W\hat{\xi}_1\perp \left[\left\{A_2,\rho_2(\Gamma) \right\}^{\prime\prime}
  P_\pm^{(2)}\mathcal{H}_2\right]=\widetilde{\mathcal{H}}_2\;\;\;\text{ and }\\
  \left(aW\hat{\xi}_1, W\hat{\xi}_1 \right)=\left(a\hat{\xi}_2,
  \hat{\xi}_2 \right) \text{ for all }\;\;
  a\in\left\{A_2,\rho_2(\Gamma)
   \right\}^{\prime\prime}.
   \end{split}
\end{eqnarray}
Define  linear map $K$ by $K\left(v \right)=
\left\{\begin{array}{ll}
a\hat{\xi}_2 ,&\text{ if }\;v=aW\hat{\xi}_1 \;\;\; a\in a\in\left\{A_2,\rho_2(\Gamma)
   \right\}^{\prime\prime},\\
0 ,&\text{ if }\;v\in \mathcal{H}_2\ominus\left[\left\{A_2,\rho_2(\Gamma)
\right\}^{\prime\prime}
  \hat{\xi}_2\right] .\end{array}\right.$
  It follows from (\ref{0.71}) that $K$ extends to the partial isometry from
$\left\{A_2,\rho_2(\Gamma)
\right\}^{\prime}$. Therefore, there exists unitary $\widetilde{K}\in\left\{A_2,\rho_2(\Gamma)
\right\}^{\prime}$ with the property: $\widetilde{K}v=Kv$ for all
$v\in\left[\left\{A_2,\rho_2(\Gamma)
\right\}^{\prime\prime}
  W\hat{\xi}_1\right]$. Thus $\mathcal{U}=\widetilde{K}W$ satisfies
  the conditions of  proposition \ref{Prop12}.
\end{proof}

\paragraph{The parameters of the states from paragraph \ref{parnatexmmpl}.}
Here we follow the notation of paragraphs  \ref{parnatexmmpl} and
\ref{paragraph2.1}.
 \subparagraph{State $\varphi_{sp}$.} Below we find  parameters
 $\left(\mathcal{H}, A,\widetilde{\mathcal{H}},\rho \right)$ from
 paragraph \ref{paragraph2.1}  such that
 $\varphi_{sp}=\psi_A^\rho$, where $\psi_A^\rho$ defined in
 proposition \ref{Prop11a}.

Let $\left(\rho,\mathcal{H}_\varphi, \xi_\varphi \right)$ be
GNS-representation of group $\Gamma$ corresponding to $\varphi$,
where $\varphi(\gamma)=\left(\rho(\gamma)\xi_\varphi,\xi_\varphi
\right)$ for all $\gamma \in\Gamma $  and $\mathcal{H}_\varphi=
\left[\rho\left(\Gamma  \right)\xi_\varphi\right]$.
  An easy computation shows that $\mathcal{H}=\mathcal{H}_\varphi$,
   $A$ acts by
  \begin{eqnarray}
A\xi= \left(\xi,\xi_\varphi \right)\xi_\varphi \;\;\;
(\xi\in\mathcal{H}),
  \end{eqnarray}
and $\widetilde{\mathcal{H}}=\mathcal{H}$. It is clear
$\mathcal{H}_{reg}=0$.
  \subparagraph{State $\varphi_{reg}$.} As above  $\left(\rho_\varphi,\mathcal{H}_\varphi, \xi_\varphi
  \right)$ is GNS-representation of  $\Gamma$. If   $\left(\rho_\varphi^{(k)},
  \mathcal{H}_\varphi^{(k)}, \xi_\varphi^{(k)}\right)$ is $k$-th
  copy of $\left(\rho_\varphi,\mathcal{H}_\varphi, \xi_\varphi
  \right)$ then $$\mathcal{H}=\mathcal{H}_{reg}=
  \bigoplus\limits_{k=1}^\infty\left(\rho_\varphi^{(k)},
  \mathcal{H}_\varphi^{(k)}, \xi_\varphi^{(k)}\right).$$
  It is obvious, $A\equiv 0$. Now define $\mathfrak{e}_{kl}^\prime$
  by
  $$\mathfrak{e}_{kl}^\prime\left(\xi_1, \xi_2,\ldots \right)=
  \left(\underbrace{0,\ldots,0}_{k-1},\xi_l,0,0,\ldots\right).$$
Put $\rho=\bigoplus\limits_{k=1}^\infty\rho_\varphi^{(k)}$,
$\hat{\xi}=\left(\xi_\varphi,0,0,\ldots \right)$. It is easy to
check that $\varphi_{reg}=\psi_0^\rho$.
\paragraph{$\mathfrak{S}_\infty$-invariance of $\psi_A^\rho$.} The
next assertion follows from definition of $\psi_A^\rho$.
\begin{Prop}\label{Prop13a}
Let $s\in\mathfrak{S}_\infty$,
$\gamma=\left(\gamma_1,\gamma_2,\ldots\right)\in\Gamma^\infty_0$. If
$s \gamma = \prod\limits_{p\in \mathbb{N}\diagup s}
 s_p \gamma(p)$, where $s_p \gamma(p)$ is generalized cycle of
 $s\gamma$ (see (\ref{product})), then
 $\psi_A^\rho\left(s\gamma \right)=\prod\limits_{p\in \mathbb{N}\diagup s}
  \psi_A^\rho\left(s_p \gamma(p) \right)$. In particular, it follows
  from Proposition \ref{multiplicativity} that $\psi_A^\rho$ is
  indecomposable state on $\Gamma\wr\mathfrak{S}_\infty$.
\end{Prop}
Denote by $\left(n_1\;\,n_2\;\,\ldots \;\, n_k \right)$ cycle $
\left\{n_1\mapsto n_2\mapsto\ldots\mapsto n_k\mapsto n_1
\right\}\in\mathfrak{S}_\infty$. Suppose that
$\gamma=\left(\gamma_1,\gamma_2,\ldots\right)\in \Gamma^\infty_e$
satisfies the condition: $\gamma_i=e$ for all $i\notin
\left\{n_1,n_2,\ldots,n_k \right\}$. If ${\rm Tr}\left(|A|
\right)=1$, $c_k=\left(n_1\;\,n_2\;\,\ldots \;\, n_k \right)$ then,
using (\ref{psik}), we have
\begin{eqnarray}\label{formula1}
\psi_A^\rho\left(c_k\gamma\right)={\rm Tr}^{\otimes N}\big(
 U\left(c_k\right)\left(\rho\left(\gamma_1\right)\otimes\rho\left(\gamma_2\right)\otimes\ldots
\otimes\rho\left(\gamma_N\right) \right)
 A^{\otimes N}\big)
\end{eqnarray}
for all  $N\geq {\rm max}\left\{n_1,n_2,\ldots,n_k \right\}$, where
${\rm Tr}^{\otimes N}$ is the ordinary trace on
$\mathcal{B}\left(\mathcal{H} \right)^{\otimes N}$, $ A^{\otimes
N}=\underbrace{A\otimes\ldots\otimes A}_N$. The next lemma extends
 formula {\ref{formula1}} on the general case.
 \begin{Lm}\label{Lm14}
If $k>1$ then
\begin{eqnarray*}
\psi_A^\rho\left(c_k\gamma\right)=
 {\rm Tr}^{\otimes N}\big(
 U\left(\left(n_1\;\,n_2\;\,\ldots \;\, n_k \right)\right)\left(\rho\left(\gamma_{n_1}\right)
 \otimes\rho\left(\gamma_{n_2}\right)\otimes\ldots
\otimes\rho\left(\gamma_{n_k}\right) \right)
 A^{\otimes k}\big).
\end{eqnarray*}
 \end{Lm}
\begin{proof}
Let $\widetilde{P}$ be an orthogonal projection on subspace
$\widetilde{\mathcal{H}}=\mathcal{H}_+\oplus\mathcal{H}_-$ (see
paragraph \ref{paragraph2.1}). Put $E=E_1\otimes
E_2\otimes\ldots\otimes E_N\otimes\ldots$, where $E_i=
\left\{\begin{array}{ll}\widetilde{P}+\mathfrak{e}^\prime_{ii},
&\textit{ if }\;i=n_j,\\
      I_\mathcal{H},&\textit{ if } i\neq n_j \text{ for all } j\in  \left\{1,2,\ldots,k \right\}.
\end{array}\right.$
Considering identical operator $I\in\mathcal{B}\left(\mathcal{H}
\right)$ as element of  $\mathcal{H}_A^\rho$, we obtain from
(\ref{psik}), (\ref{psik1}), (\ref{psik2})
\begin{eqnarray}\label{EII}
E I=I.
\end{eqnarray}
It follows from (\ref{UmatrixUnit}) that
\begin{eqnarray}\label{UEU}
\widetilde{E}=U\left(c_k \right)EU\left(c_k \right)^*E
=\widetilde{E}_1\otimes
\widetilde{E}_2\otimes\ldots\otimes\widetilde{E}_N\otimes\ldots,
\end{eqnarray}
where $\widetilde{E}_i=\left\{\begin{array}{ll}\widetilde{P},
&\textit{ if }\;i=n_j,\\
      I_\mathcal{H},&\textit{ if } i\neq n_j \text{ for all } j\in  \left\{1,2,\ldots,k \right\}.
\end{array}\right.$
By properties (1)-(4) from paragraph \ref{paragraph2.1}, using
(\ref{Piarho}) and (\ref{UmatrixUnit}), we obtain
\begin{eqnarray}\label{relee}
\Pi_A^\rho\left(\gamma\right)E=E\Pi_A^\rho\left(\gamma\right),\;
\Pi_A^\rho\left(\gamma\right)\widetilde{E}
=\widetilde{E}\Pi_A^\rho\left(\gamma\right).
\end{eqnarray}
Thus
\begin{eqnarray}
\begin{split}
 \psi_A^\rho\left(c_k\gamma\right)=
 \left(\Pi_A^\rho\left(c_k \right)\Pi_A^\rho(\gamma )I,I
 \right)\stackrel{(\ref{EII})}{=}
 \left(\Pi_A^\rho\left(c_k \right)\Pi_A^\rho(\gamma )EI,EI \right)\\
 =\left(\Pi_A^\rho\left(c_k \right)\Pi_A^\rho(\gamma )
 \Pi_A^\rho\left(c_k \right)^*\left[\Pi_A^\rho\left(c_k \right)E
  \Pi_A^\rho\left(c_k \right)^*\right]\Pi_A^\rho\left(c_k \right)I,EI
  \right)\\
  \stackrel{(\ref{relee})}{=}
  \left(\Pi_A^\rho\left(c_k \right)\Pi_A^\rho(\gamma )
 \Pi_A^\rho\left(c_k \right)^*\Pi_A^\rho\left(c_k \right)I,\left[\Pi_A^\rho\left(c_k \right)E
  \Pi_A^\rho\left(c_k \right)^*\right]EI \right)\\
  \stackrel{(\ref{UEU})}{=}
  \left(\Pi_A^\rho\left(c_k \right)\Pi_A^\rho(\gamma )I,
 \widetilde{E}I \right)\stackrel{(\ref{UEU}),(\ref{UmatrixUnit})}{=}
  \left(\Pi_A^\rho\left(c_k \right)\Pi_A^\rho(\gamma )\widetilde{E}I,
 \widetilde{E}I \right).
 \end{split}
\end{eqnarray}
Hence, applying (\ref{psik}), (\ref{psik1}), (\ref{psik2}), obtain
for $N\geq {\rm max}\left\{n_1,n_2,\ldots,n_k \right\}$
 $\psi_A^\rho\left(c_k\gamma\right)=\,_1\psi_N
 \left(\widetilde{E}U\left(c_k\right)
 \left(\rho\left(\gamma_1\right)\otimes\rho
 \left(\gamma_2\right)\otimes\ldots
\otimes\rho\left(\gamma_N\right) \right)\widetilde{E}\right)$. Since
$\widetilde{P}\perp\mathfrak{e}_{kk}^\prime$ for all $k$, then
$_1\psi_N\left(\widetilde{E}U\left(c_k\right)
 \left(\rho\left(\gamma_1\right)\otimes\rho\left(\gamma_2\right)\otimes\ldots
\otimes\rho\left(\gamma_N\right) \right)\widetilde{E}\right)\\=
 {\rm Tr}^{\otimes N}\big(
 U\left(\left(n_1\;\,n_2\;\,\ldots \;\, n_k \right)\right)\left(\rho\left(\gamma_{n_1}\right)
 \otimes\rho\left(\gamma_{n_2}\right)\otimes\ldots
\otimes\rho\left(\gamma_{n_k}\right) \right)
 A^{\otimes k}\big)$.
\end{proof}
\begin{Rem}
One should notice that in the case in which $c_k=1$,
\begin{eqnarray}
\psi_A^\rho\left(\gamma  \right)=\prod\limits_{n=1}^\infty
 \left[{\rm Tr}\left(\rho\left(\gamma_n  \right)|A| \right)+\left(1-{\rm
Tr}\left(|A| \right) \right)\left(\rho\left(\gamma_n
\right)\hat{\xi},\hat{\xi} \right) \right].
\end{eqnarray}
\end{Rem}
Hence, taking into account Proposition \ref{Prop13a}, Lemma
\ref{Lm14} and (\ref{formula1}), we obtain the next important
property
\begin{eqnarray}\label{sinfinv}
\psi_A^\rho\left(sgs^{-1} \right)=\psi_A^\rho\left(g \right) \text{
for all } s\in\mathfrak{S}_\infty, g\in \Gamma\wr
\mathfrak{S}_\infty.
\end{eqnarray}
\subsubsection{KMS-condition for the $\mathfrak{S}_\infty$-central states.}
\label{KMSsec}
\paragraph{KMS-condition for $\psi_A^\rho$.}
To the general definition of the KMS-condition we refer the reader
to the book \cite{Tak}. Here we introduce the definition of the
KMS-condition for the indecomposable states only.

\begin{Def}\label{KMS} Let $\varphi$ be an indecomposable state on the group
$G$. Let $\left(\pi_\varphi,\mathcal{H}_\varphi,\xi_\varphi\right)$ be the
corresponding GNS-construction, where $\xi_\varphi$ is such that
$\varphi(g)=\left(\pi_\varphi(g)\xi_\varphi,\xi_\varphi\right)$ for each
$g\in G$. We say that $\varphi$ satisfies  the KMS-condition or $\varphi$
is KMS-state, if $\xi_\varphi$ is separating\footnote{This means that for
every $a\in\pi_\varphi(G)^{\prime\prime}$ the conditions $a\xi_\varphi =0$
and $a=0$ are equivalent.} for the $w^*$-algebra
$\pi_\varphi(G)^{\prime\prime}$, generated by operators $\pi_\varphi(G)$.
\end{Def}
The main result of this paragraph is the following:
\begin{Th}\label{theorem15}
Let $\left(A, \hat{\xi}, \mathcal{H}_{reg}, \mathfrak{e}_{kl}^\prime
\right)$ satisfy the conditions {\rm (1)}-{\rm (4)} from paragraph
\ref{paragraph2.1}. State $\psi_A^\rho$ satisfies the KMS-condition
if and only if ${\rm Ker}\, A=\mathcal{H}_{reg}$ and $\hat{\xi}$ is
cyclic and separating for the restriction
$\rho_{11}=\rho\Big|_{\mathfrak{e}_{11}^\prime\mathcal{H}_{reg}}$ of
representation $\rho$ to subspace
$\mathfrak{e}_{11}^\prime\mathcal{H}$.
\end{Th}

As a preliminary to the proof of the theorem, we will discuss two
auxiliary lemmas.
\begin{Lm}\label{lemma16}
Let $\left(\pi_{\psi_k},H_{\psi_k},\xi_{\psi_k} \right)$ be
GNS-representation of $\mathcal{B}\left(\mathcal{H} \right)$ corresponding
to state $\psi_k$ (see (\ref{psik})). Fix any $\epsilon >0$ and denote by
$P_{[\epsilon,1 ]}$ the spectral projection of $\left|A \right|$. Then for
each $a\in\mathcal{B}\left(\mathcal{H} \right)$ the map
\begin{eqnarray*}
\mathfrak{R}_ {P_{[\epsilon,1 ]}aP_{[\epsilon,1
]}}:x\mapsto x\cdot P_{[\epsilon,1
]}\,a\,P_{[\epsilon,1 ]}
\end{eqnarray*}
 may be extended
by continuous to the bounded operator on $H_{\psi_k}$ and
 $\left\| \mathfrak{R}_ {P_{[\epsilon,1 ]}aP_{[\epsilon,1
]}} \right\|_{H_{\psi_k}}\leq \frac{\left\|a
\right\|}{\sqrt{\epsilon} }$.
\end{Lm}
\begin{proof}
Put $b=P_{[\epsilon,1 ]}aP_{[\epsilon,1 ]}$. Then
\begin{eqnarray*}
&\left(\mathfrak{R}_bx,\mathfrak{R}_bx \right)_{H_{\psi_k}}= {\rm
Tr}\left(b|A|b^*x^*x \right)\leq \left\| b|A|b^*\right\|  {\rm
Tr}\left( P_{[\epsilon,1 ]}x^*x\right)\\
&=\left\| b|A|b^*\right\|\cdot{\rm Tr}\left(|A|
\cdot\left[\sum\limits_{\lambda\in[\epsilon,1 ]\cap\,{\rm
Spectrum}\,|A|} \lambda^{-1}P_\lambda\right]x^*x\right)\\
&\leq\epsilon^{-1}\cdot\left\| b|A|b^*\right\|\cdot{\rm
Tr}\left(|A|P_{[\epsilon,1 ]}x^*x
\right)\leq\epsilon^{-1}\cdot\left\| b|A|b^*\right\|\cdot{\rm
Tr}\left(|A|x^*x\right)\leq\\
&\stackrel{\ref{psik}}{=}\epsilon^{-1}\cdot\left\|
b|A|b^*\right\|\psi_k\left( x^*x\right)\leq\epsilon^{-1}\cdot\left\|
b\right\|^2\left(x^*x \right)_{H_{\psi_k}}.
\end{eqnarray*}
\end{proof}
\begin{Lm}\label{lemma17}
Suppose that for $\left(A, \hat{\xi}, \mathcal{H}_{reg},
\mathfrak{e}_{kl}^\prime \right)$  the conditions {\rm (1)}-{\rm (4)} from
paragraph \ref{paragraph2.1} hold. Denote by $P_0$ and $P_{reg}$ the
orthogonal projections onto ${\rm Ker}\,A$ and $\mathcal{H}_{reg}$
respectively. Let $ \left[\Pi_A^\rho\left(\Gamma\wr\mathfrak{S}_\infty
\right) I\right]$ be the subspace in $\mathcal{H}_A^\rho$ (see paragraphs
\ref{paragraph2.2}, \ref{paragraph2.3}), generated by
$\Pi_A^\rho\left(\Gamma\wr\mathfrak{S}_\infty \right) I$. For $m\in
 \left\{\rho\left(\Gamma \right) \right\}^\prime\subset\mathcal{B}(\mathcal{H})$
define the linear map
$\mathfrak{R}_m^{(k)}:\mathcal{B}(\mathcal{H})^{\otimes\infty}\mapsto\mathcal{B}(\mathcal{H})^{\otimes\infty}$
as follows
\begin{eqnarray}
\begin{split}
\mathfrak{R}_m^{(k)}\left(a_1\otimes\ldots\otimes a_k\otimes
a_{k+1}\otimes\ldots \right)\\
=a_1\otimes\ldots\otimes a_k\cdot
\mathfrak{e}_{kk}^\prime\cdot m\cdot \mathfrak{e}_{kk}^\prime\otimes a_{k+1}\otimes\ldots.
\end{split}
\end{eqnarray}
If $P_0=P_{reg}$ then
\begin{itemize}
    \item {\rm (i)} $\mathfrak{R}_m^{(k)}\left(
        \Pi_A^\rho\left(\Gamma\wr\mathfrak{S}_\infty \right)
        I\right)\subset\left[\Pi_A^\rho\left(\Gamma\wr\mathfrak{S}_\infty
        \right) I\right]$;
    \item {\rm (ii)} the extension of $\mathfrak{R}_m^{(k)}
        \Big|_{\Pi_A^\rho\left(\Gamma\wr\mathfrak{S}_\infty \right)I}$
        by continuous is bounded operator in
        $\left[\Pi_A^\rho\left(\Gamma\wr\mathfrak{S}_\infty \right)
        I\right] \subset\mathcal{H}_A^\rho$.
\end{itemize}
\end{Lm}
\begin{proof}
To prove {\rm (i)}, it suffices to show that
$\mathfrak{R}_m^{(k)}(I)\in
 \left[\Pi_A^\rho\left(\Gamma\wr\mathfrak{S}_\infty
        \right) I\right]$. Indeed, by property {\rm (4)},
for any $\epsilon
>0$ there exists $a_\epsilon =
 \sum\limits_{g\in \Gamma_\epsilon } c_\gamma \rho(\gamma)$, where $\Gamma_\epsilon$ is
a finite subset in $\Gamma$, satisfying
\begin{eqnarray*}
\left\| \mathfrak{e}_{1k}^\prime m \mathfrak{e}_{k1}^\prime\hat{\xi}-
 a_\epsilon \hat{\xi}\right\|_\mathcal{H}<\epsilon.
\end{eqnarray*}
Hence, considering $\mathfrak{R}_m^{(k)}(I)$ and $a_\epsilon^{(k)}=
 \underbrace{I\otimes\ldots\otimes I}_{k-1}\otimes P_{reg}a_\epsilon P_{reg}\otimes I\otimes\ldots $
 as the elements from $\mathcal{H}_A^\rho$, we have
\begin{eqnarray}\label{rightm}
\left\| \mathfrak{R}_m^{(k)}(I)- a_\epsilon^{(k)}  \right\|_{\mathcal{H}_A^\rho}<\epsilon.
\end{eqnarray}
It follows from (\ref{astrans}) and (\ref{Oaction}), that operator of the
left multiplication on $\underbrace{I\otimes\ldots\otimes I}_{k-1}\otimes
P_0\otimes I\otimes\ldots $ lies in
 $\Pi_A^\rho\left(\Gamma\wr\mathfrak{S}_\infty \right)^{\prime\prime}$.
Hence, since $P_0=P_{reg}$, we get  $a_\epsilon^{(k)}\in
\Pi_A^\rho\left(\Gamma\wr\mathfrak{S}_\infty \right)^{\prime\prime}$.
Therefore, using (\ref{rightm}),  we obtain $\mathfrak{R}_m^{(k)}(I)\in
 \left[\Pi_A^\rho\left(\Gamma\wr\mathfrak{S}_\infty
        \right) I\right]$.

        Let us prove statement {\rm (ii)}.
        Put $\mathfrak{S}_\infty^{(k)}=
        \left\{s\in\mathfrak{S}_\infty:s(k)=k \right\}$.
        First, using (\ref{sinfinv}),
         we observe that
         \begin{eqnarray}\label{centraliser}\begin{split}
\left(a_1b_1I,a_2b_2I\right)_{\mathcal{H}_A^\rho}=
\left(a_1b_1b_2^*0I,a_2I\right)_{\mathcal{H}_A^\rho}\\
\text{ for all
}\; a_1, a_2\in\Pi_A^\rho\left(\Gamma\wr\mathfrak{S}_\infty
        \right)^{\prime\prime} \text{ and }\;
        b_1, b_2\in\Pi_A^\rho\left(\mathfrak{S}_\infty
        \right)^{\prime\prime}.
       \end{split}  \end{eqnarray}
 Denote be $\mathcal{L}_{P_0}^{(k)}$ operator of the
left multiplication on $\underbrace{I\otimes\ldots\otimes I}_{k-1}\otimes
P_0\otimes I\otimes\ldots $. By
 (\ref{astrans}) and (\ref{Oaction}),  $\mathcal{L}_{P_0}^{(k)}\in
 \Pi_A^\rho\left(\mathfrak{S}_\infty\right)^{\prime\prime}$.
 Therefore, $ \left[\Pi_A^\rho\left(\Gamma\wr\mathfrak{S}_\infty
        \right)\left(I- \mathcal{L}_{P_0}^{(k)}\right)I
        \right]$, $\mathbf{H}_l= \left[\Pi_A^\rho
        \left(\left(k\;\,l \right)\cdot\mathfrak{S}_\infty^{(k)} \right)
        \Pi_A^\rho\left(\Gamma^\infty_e\right)
        \mathcal{L}_{P_0}^{(k)}I \right]$
        $\left(l\in\mathbb{N} \right)$ are the subspaces in
        $ \left[\Pi_A^\rho\left(\Gamma\wr\mathfrak{S}_\infty
        \right)I
        \right]$ and, according to (\ref{centraliser}), we have
\begin{eqnarray}\label{perp1}
\left[\Pi_A^\rho\left(\Gamma\wr\mathfrak{S}_\infty
        \right)\left(I- \mathcal{L}_{P_0}^{(k)}\right)I
        \right] \perp \mathbf{H}_l \text{ for all }
        l\in\mathbb{N}.
\end{eqnarray}
Now we prove that subspaces $ \left\{\mathbf{H}_l
\right\}_{l\in\mathbb{N}}$ are pairwise orthogonal. For convenience
we assume that $k=1$. Denote by $E_m$ the orthogonal projection on
subspace
$\mathbb{C}\mathfrak{e}_{m1}^\prime\hat{\xi}\subset\mathcal{H}$
$\left(m\in\mathbb{N} \right)$. Put $A_m=A+\left(I-{\rm Tr}\,|A|
\right)E_m$,
$\mathfrak{E}_m^{(i)\prime}=\underbrace{I\otimes\ldots\otimes
I}_{i-1}\otimes \mathfrak{e}_{mm}^\prime\otimes I\otimes\ldots$ and
$E_m^{(i)}=\underbrace{I\otimes\ldots\otimes I}_{i-1}\otimes
E_m\otimes I\otimes\ldots$. By definition,
\begin{eqnarray}\label{ee}
E_m^{(i)}\mathfrak{E}_l^{(i)\prime}=\delta_{ml}E_m^{(i)}, \text{
where }\delta_{ml} \text{ is Kronecker's delta}.
\end{eqnarray}
  It follows from   the definition of $A_m$ that for $s^{-1}(1)\neq
1$ and $n>s^{-1}(1)$
\begin{eqnarray}\label{zero}
 \mathfrak{E}_1^{\left(s^{-1}(1)\right)\prime}\cdot\bigotimes_{m=1}^n
A_m=0.
\end{eqnarray}
 Fix any
$\widetilde{\gamma }, \widehat{\gamma }\in\Gamma^\infty_e$, $ s_1\in
\left(1\,\;l_1 \right)\mathfrak{S}_\infty^{(1)}$ and $s_2\in
\left(1\,\;l_2 \right)\mathfrak{S}_\infty^{(1)}$. Let us show  that
for $l_1\neq l_2$
\begin{eqnarray}
\kappa=\left(\Pi_A^\rho\left(s_1\widetilde{\gamma }
\right)\mathcal{L}_{P_0}^{(1)}I,
 \Pi_A^\rho\left(s_2\widehat{\gamma }
\right)\mathcal{L}_{P_0}^{(1)}I
\right)_{\mathcal{H}_A^\rho}=0.
\end{eqnarray}
Let ${\rm Tr}^{\otimes n}$ be the ordinary trace on $w^*$-factor
$\mathcal{B}\left(\mathcal{H} \right)^{\otimes n}$.
 If $s=s_2^{-1}s_1$,
$\gamma_m=\widehat{\gamma}_{s(m)}^{\,-1}\cdot\widetilde{\gamma}_m\in\Gamma
$, $\gamma=\left(\gamma_1, \gamma_2,\ldots\right)$ and $n> {\rm max}
\left\{{\rm max}\left\{i:\gamma_i\neq e  \right\},{\rm max}
\left\{i:s(i)\neq i \right\}\right\}$  then, using definition of
$\Pi_A^\rho$ (see (\ref{Piarho})), we have
\begin{eqnarray}
\kappa={\rm Tr}^{\otimes n}\left(E_1^{(1)}\cdot
U_n(s)\cdot\bigotimes_{m=1}^n\rho \left(\gamma_m \right) \cdot
E_1^{(1)}\cdot\bigotimes_{m=1}^n A_m\right),
\end{eqnarray}
where $U_n(s)$ is defined in paragraph \ref{paragraph2.3}. Hence,
applying property {\rm (4)} from paragraph {\ref{paragraph2.1}},
(\ref{ee}) and (\ref{UmatrixUnit}), we obtain
\begin{eqnarray*}
\kappa={\rm Tr}^{\otimes n}\left(E_1^{(1)}\cdot
U_n(s)\,\left(U_n(s)\right)^*\mathfrak{E}_1^{(1)\prime}\,
U_n(s)\cdot\bigotimes_{m=1}^n\rho \left(\gamma_m \right) \cdot
E_1^{(1)}\cdot\bigotimes_{m=1}^n A_m\right)\\
\stackrel{(\ref{UmatrixUnit})}{=}{\rm Tr}^{\otimes
n}\left(E_1^{(1)}\cdot
U_n(s)\,\mathfrak{E}_1^{\left(s^{-1}(1)\right)\prime}\,
\cdot\bigotimes_{m=1}^n\rho \left(\gamma_m \right) \cdot
E_1^{(1)}\cdot\bigotimes_{m=1}^n A_m\right)\\
\stackrel{{\rm property (4)}}{=}{\rm Tr}^{\otimes
n}\left(E_1^{(1)}\cdot U_n(s) \cdot\bigotimes_{m=1}^n\rho
\left(\gamma_m \right) \cdot
E_1^{(1)}\cdot\mathfrak{E}_1^{\left(s^{-1}(1)\right)\prime}\cdot\bigotimes_{m=1}^n
A_m\right)\stackrel{(\ref{zero})}{=} 0.
\end{eqnarray*}
Therefore,
\begin{eqnarray}\label{perp2}
\mathbf{H}_l\perp \mathbf{H}_m \text{ for all } l\neq m.
\end{eqnarray}
As in the proof of {\rm (i)},
$\mathfrak{R}_m^{(1)}(I)=\mathfrak{e}_{11}^\prime m\mathfrak{e}_{11}^\prime
\otimes I\otimes I\otimes\ldots$ lies in subspace $
\left[\Pi_A^\rho\left(\Gamma_e^\infty \right)
\mathfrak{L}_{P_0}^{(1)}I\right]\subset \mathbf{H}_1$. Therefore,
\begin{eqnarray}
\Pi_A^\rho
        \left(\left(1\;\,l \right)\cdot\mathfrak{S}_\infty^{(1)} \right)
        \Pi_A^\rho\left(\Gamma^\infty_e\right)
        \mathcal{L}_{P_0}^{(1)}\mathfrak{R}_m^{(1)}(I)\subset\mathbf{H}_l.
\end{eqnarray}
Further, using (\ref{UmatrixUnit}) and relation
\begin{eqnarray*}
\mathfrak{R}_m^{(1)}\Pi_A^\rho
        \left(\left(1\;\,l \right)\cdot s \right)
        \Pi_A^\rho(\gamma )
        \mathcal{L}_{P_0}^{(1)}
        (I)\stackrel{(\ref{UmatrixUnit})}{=}
 \mathcal{L}_{\mathfrak{e}_{11}^\prime m\mathfrak{e}_{11}^\prime}^{(l)}
 \Pi_A^\rho
        \left(\left(1\;\,l \right)\cdot s \right)
        \Pi_A^\rho(\gamma )
        \mathcal{L}_{P_0}^{(1)}(I),
\end{eqnarray*}
where $s\in\mathfrak{S}_\infty^{(1)}$, $\gamma \in\Gamma^\infty_e$,
we obtain that $\mathfrak{R}_m^{(1)}$ is the bounded operator on
$\mathbf{H}_l$ and $\left\|\mathfrak{R}_m^{(1)}
\right\|_{\mathbf{H}_1}\leq \left\|\mathfrak{e}_{11}^\prime
m\mathfrak{e}_{11}^\prime \right\|_{\mathcal{H}}$. Since, by
 (\ref{perp1}) and (\ref{perp2}),
 \begin{eqnarray}
\left[\Pi_A^\rho\left(\Gamma\wr\mathfrak{S}_\infty
        \right)I
        \right]=
 \left[\Pi_A^\rho\left(\Gamma\wr\mathfrak{S}_\infty
        \right)\left(I- \mathcal{L}_{P_0}^{(1)}\right)I
        \right]\bigoplus\limits_{m=1}^\infty \mathbf{H}_m,
 \end{eqnarray}
 and $\left[\Pi_A^\rho\left(\Gamma\wr\mathfrak{S}_\infty
        \right)\left(I- \mathcal{L}_{P_0}^{(1)}\right)I
        \right]\subset {\rm Ker}\,\mathfrak{R}_m^{(1)}$, operator
        $\mathfrak{R}_m^{(1)}$ is bounded on subspace $\left[\Pi_A^\rho\left(\Gamma\wr\mathfrak{S}_\infty
        \right)I
        \right]$.
\end{proof}
\begin{proof}[{\bf The proof of Theorem \ref{theorem15}.}]
Let $\Pi_A^{\rho\,0}$ be the restriction $\Pi_A^\rho$ to subspace $
\left[ \Pi_A^\rho\left(\Gamma\wr\mathfrak{S}_\infty
\right)I\right]$. Obvious, $\Pi_A^{\rho\,0}$ and GNS-representation
of $\Gamma\wr\mathfrak{S}_\infty$, corresponding to $\psi_A^\rho$,
are naturally unitary equivalent. Let us prove that $I$ is the
cyclic vector for $\Pi_A^{\rho\,0}\left(\Gamma\wr\mathfrak{S}_\infty
\right)^\prime$.

For any $n\in\mathbb{N}$ fix $\gamma
=\left(\gamma_1,\gamma_2,\ldots,\gamma_n,e,e,\ldots\right)\in\Gamma^\infty_e$
and $s\in\mathfrak{S}_n$. Put $\eta=\Pi_A^\rho(\gamma)I=
 \left(\bigotimes\limits_{m=1}^n \rho\left(\gamma_m\right)\right)\otimes I\otimes
 I\otimes\ldots\in  \left[\Pi_A^{\rho\,0}\left(\Gamma^\infty_e
\right)I \right]\subset
\left[\Pi_A^{\rho\,0}\left(\Gamma\wr\mathfrak{S}_\infty \right)I \right]$.
If $P_{[\epsilon,1 ]}$ is the spectral projection of $|A|$ then, by
(\ref{astrans}), (\ref{Oaction}) and lemma \ref{lemma17} (i), for every
$m_j^\prime\in \rho(\Gamma)^\prime$
\begin{eqnarray*}
a_\epsilon
=\left(\bigotimes\limits_{j=1}^n\left(P_{[\epsilon,1
]}\rho\left(\gamma_j\right)P_{[\epsilon,1
]}+\mathfrak{e}_{jj}^\prime m_j^\prime\mathfrak{e}_{jj}^\prime
\right) \right)\otimes I\otimes  I\otimes
\ldots\in\left[\Pi_A^{\rho\,0}\left(\Gamma\wr\mathfrak{S}_\infty
\right)I \right].
\end{eqnarray*}
Since $\hat{\xi}$ is cyclic and separating for the restriction
$\rho_{11}=\rho\Big|_{\mathfrak{e}_{11}^\prime\mathcal{H}_{reg}}$
and ${\rm Ker}\, A=\mathcal{H}_{reg}$, then for any $\delta >0$
 there exist $\epsilon >0$ and
  $ \left\{m_j^\prime \right\}_{j=1}^n\subset\rho(\Gamma)^\prime$
  such that
  \begin{eqnarray*}
\left\| \Pi_A^\rho(\gamma) I  - a_\epsilon
\right\|_{\mathcal{H}_A^\rho}<\delta.
  \end{eqnarray*}
  But, by lemmas \ref{lemma16}-\ref{lemma17}, operator
  $\mathfrak{R}_{a_\epsilon }$ of right multiplication on
  $a_\epsilon$ lies in $\Pi_A^{\rho\,0}\left(\Gamma\wr\mathfrak{S}_\infty
\right)^\prime$. Therefore,
\begin{eqnarray}\label{gammaprim}
\Pi_A^\rho(\gamma) I\in
\left[\Pi_A^{\rho\,0}\left(\Gamma\wr\mathfrak{S}_\infty
\right)^\prime I \right].
\end{eqnarray}
Now we note that, by (\ref{sinfinv}), the right multiplication on
$U(s)$ defines the unitary operator
$\mathfrak{R}_{U(s)}\in\Pi_A^{\rho\,0}\left(\Gamma\wr\mathfrak{S}_\infty
\right)^\prime$. It follows from (\ref{gammaprim}) that
$\Pi_A^\rho\left(\gamma\,s \right)I=
\mathfrak{R}_{U(s)}\left(\Pi_A^\rho(\gamma) I
\right)\in\left[\Pi_A^{\rho\,0}\left(\Gamma\wr\mathfrak{S}_\infty
\right)^\prime I \right]$. Therefore $I$ is the cyclic vector for
$\Pi_A^{\rho\,0}\left(\Gamma\wr\mathfrak{S}_\infty \right)^\prime$.

Conversely, suppose that $\psi_A^\rho$ is KMS-state on
$\Gamma\wr\mathfrak{S}_\infty$. Define state $\widehat{\psi}_A^\rho\in
 \Pi_A^{\rho\,0}\left(\Gamma\wr\mathfrak{S}_\infty \right)^{\prime\prime}_*$
as follows
\begin{eqnarray}
\widehat{\psi}_A^\rho(a)=\left(a I,I  \right)_{\mathcal{H}_A^\rho}.
\end{eqnarray}
Then, by propositions \ref{multiplicativity} and \ref{Prop13a},
$\widehat{\psi}_A^\rho$ is faithful state. This means that for every $a\in
\Pi_A^{\rho\,0}\left(\Gamma\wr\mathfrak{S}_\infty \right)^{\prime\prime}$
the conditions $\widehat{\psi}_A^\rho(a^*a) =0$ and $a=0$ are equivalent.

Let us prove that ${\rm Ker}\,A=\mathcal{H}_{reg}$. If
$\mathcal{H}_{reg}\varsubsetneqq {\rm Ker}\,A$ then, by properties {\rm
(1)}-{\rm (4)} from paragraph \ref{paragraph2.1}, there exists $\gamma
\in\Gamma$ such that
\begin{eqnarray}
 \rho(\gamma)\left(P_{]0,1]}+P_{[-1,0[}\right)\neq
 \left(P_{]0,1]}+P_{[-1,0[}\right)\rho(\gamma).
\end{eqnarray}
It follows from this
\begin{eqnarray*}
Q=\left(\left(P_{]0,1]}+P_{[-1,0[}\right)\vee
 \rho(\gamma)\left(P_{]0,1]}+P_{[-1,0[}\right)\rho(\gamma )^*\right)-
 \left(P_{]0,1]}+P_{[-1,0[}\right)\neq 0.
\end{eqnarray*}
Since $Q\in \mathfrak{A}$, where $\mathfrak{A}$ is defined in property {\rm
(1)} from paragraph \ref{paragraph2.1}, then, by
(\ref{astrans})-(\ref{Oaction}), operator $\mathfrak{L}_Q^{(k)}$ of the
left multiplication on $\left(\otimes_{m=1}^{k-1}I\right)\otimes Q\otimes
I\otimes\ldots$ lies in $\Pi_A^{\rho\,0}\left(\Gamma\wr\mathfrak{S}_\infty
\right)^{\prime\prime}$. Thus  $\widehat{\psi}_A^\rho \left(
\mathfrak{L}_Q^{(k)} \right)= {\rm Tr}\,(Q\cdot|A|)=0$. But this
contradicts the faithfulness of $\widehat{\psi}_A^\rho$.

Now we prove that $\hat{\xi}$ is cyclic and separating for the
representation
$\rho_{11}=\rho\Big|_{\mathfrak{e}_{11}^\prime\mathcal{H}_{reg}}$.
Denote by $E_{11}$ the projection onto $\left[
\rho_{11}(\Gamma)^\prime\hat{\xi}\right]$ and suppose $\left[
\rho_{11}(\Gamma)^\prime\hat{\xi}\right]\subsetneqq\left[
\rho_{11}(\Gamma)\hat{\xi}\right]$. It follows from this that
\begin{eqnarray}\label{nokms}
E_{11}\in\rho_{11}(\Gamma)^{\prime\prime},\;
F_{11}=\mathfrak{e}_{11}^\prime - E_{11}\neq 0 \;\text{ and } \;
F_{11}\hat{\xi}=0.
\end{eqnarray}
Denote by $P_{reg}$ the orthogonal projection onto
$\mathcal{H}_{reg}$. Since ${\rm Ker}\, A=\mathcal{H}_{reg}$, then
\begin{eqnarray*}
P_{reg}\in\mathfrak{A}\; \text{ and } \;P_{reg}\cdot\rho\left(\Gamma
\right)^{\prime\prime}\cdot P_{reg}\subset\mathfrak{A}.
\end{eqnarray*}
 Hence, by properties {\rm (2)} and {\rm (4)} from paragraph
 \ref{paragraph2.1}, we obtain
 \begin{eqnarray*}
F=\sum\limits_{m=1}^\infty \mathfrak{e}_{m1}^\prime\cdot F_{11}
\cdot \mathfrak{e}_{1m}^\prime\in P_{reg}\cdot\rho\left(\Gamma
\right)^{\prime\prime}.
 \end{eqnarray*}
Hence, using (\ref{astrans})-(\ref{Oaction}), we obtain that
operator $\mathfrak{L}_F^{(k)}$ of the left multiplication on
$\left(\otimes_{m=1}^{k-1}I\right)\otimes F\otimes I\otimes\ldots$
lies in $\Pi_A^{\rho\,0}\left(\Gamma\wr\mathfrak{S}_\infty
\right)^{\prime\prime}$. It follows from this and (\ref{nokms}) that
$\widehat{\psi}_A^\rho \left( \mathfrak{L}_F^{(k)} \right)=0$.
\end{proof}

\subsubsection{The main result.}\label{Themainresult}
In this section we prove the main
result of this paper:
\begin{Th}\label{main} Let $\varphi$ be any  indecomposable $\mathfrak{S}_\infty$-central
state on the group $\Gamma\wr \mathfrak{S}_\infty$. Then there exist
self-adjoint operator $A$  of the {\it trace class} (see
{\cite{RS}}) from $\mathcal{B}(\mathcal{H})$ and unitary
representation $\rho$ with the properties {\rm (1)}-{\rm (4)}
(paragraph \ref{paragraph2.1}) such that $\varphi =\psi_A^\rho$ (see
Proposition {\ref{Prop11a}}).
\end{Th}
We have divided the proof into a sequence of lemmas and
propositions. First we introduce some new objects and notations.
\paragraph{Asymptotical transposition.} Let
$\left(\pi_\varphi,H_\varphi,\xi_\varphi\right)$ be
GNS-representation of $\Gamma\wr\mathfrak{S}_\infty$ associated with
$\varphi$, where
 $\varphi(g)=\left(\pi_\varphi (g)\xi_\varphi,\xi_\varphi   \right)$
 for all $g\in \Gamma\wr\mathfrak{S}_\infty$. In the sequel for
 convenience we denote group $\Gamma\wr\mathfrak{S}_\infty$ by $G$.
 Put
 \begin{eqnarray*}
G_n(\infty)=&\Big\{ s\gamma\in G\big|\;s\in\mathfrak{S}_\infty,
 \gamma=\left(\gamma_1,\gamma_2,\ldots\right)\in\Gamma^\infty_e,&\\
&s(l)=l \;\text{ and }\;
 \gamma_l=e \;\text{ for } \;l=1,2,\cdots, n \Big\},&\\
G_n = &\left\{s\gamma\in G\big|\;s(l)=l \;\text{ and }\;
 \gamma_l=e \;\text{ for all } \;l>n \right\},&\\
 G^{(k)}= &\left\{s\gamma\in G\big| \;s(k)=k
 \;\text{ and }\;\gamma_k=e\right\}.&
 \end{eqnarray*}
It is clear that $G_0(\infty)=G$.
\begin{Prop}\label{Prop19}
Let $\left(i\;j \right)$ denotes the transposition exchanging $i$
and $j$. In the weak operator topology there exists
$\lim\limits_{j\to\infty}\pi_\varphi \left(\left(i\;j \right)
\right)$.
\end{Prop}
\begin{proof}
It is suffices to show that for any $g,h\in G$ there exists
$\lim\limits_{j\to\infty}\left(\pi_\varphi \left(\left(i\;j
\right)\right)\pi_\varphi (g)\xi_\varphi , \pi_\varphi
(g)\xi_\varphi \right)$.
 Find $N>i$ such that $g,h\in G_n$ for all $n\geq N$. Since $\varphi$
 is $\mathfrak{S}_\infty$-central, then
 \begin{eqnarray*}
\left(\pi_\varphi \left(\left(i\;N \right)\right)\pi_\varphi
(g)\xi_\varphi , \pi_\varphi (g)\xi_\varphi \right)\\=
\left(\pi_\varphi \left(\left(i\;j \right)\right)\pi_\varphi
(g)\pi_\varphi ((n\;N))\xi_\varphi , \pi_\varphi (g)\pi_\varphi
((n\;N))\xi_\varphi \right)\\= \left(\pi_\varphi \left(\left(i\;n
\right)\right)\pi_\varphi (g)\xi_\varphi , \pi_\varphi
(g)\xi_\varphi \right).
 \end{eqnarray*}
Thus $\lim\limits_{j\to\infty}\left(\pi_\varphi \left(\left(i\;j
\right)\right)\pi_\varphi (g)\xi_\varphi , \pi_\varphi
(g)\xi_\varphi \right)=\left(\pi_\varphi \left(\left(i\;N
\right)\right)\pi_\varphi (g)\xi_\varphi , \pi_\varphi
(g)\xi_\varphi \right)$.
\end{proof}
We will call $\mathcal{O}_i=\lim\limits_{j\to\infty}\pi_\varphi
\left(\left(i\;j \right) \right)$ the {\it asymptotical}
transposition.
\paragraph{The properties of the asymptotical transposition.}
\begin{Lm}\label{Fn}
  Let $g,h\in G^{(n)}$. Then for each $k\neq n$ the next relation
  holds: \begin{eqnarray}
   \left(\pi_\varphi\left(g\cdot(n\;k)\cdot h\right)\xi_\varphi,
   \xi_\varphi\right)=\left(\pi_\varphi(g)
   \mathcal{O}_k\pi_\varphi(h)\xi_\varphi,\xi_\varphi\right)
  \end{eqnarray}
\end{Lm}
\begin{proof}
 Fix  $N\in\mathbb{N}$ such that $g,h\in G_N\cap G^{(n)}$. Then for each $m>N$ we have:
 $(n\;m)\cdot g=g\cdot(n\;m)$, $(n\;m)\cdot h=h\cdot(n\;m)$. Hence, by the
  $\mathfrak{S}_\infty$-centrality of $\varphi $, we obtain
  \begin{eqnarray*}
  \left(\pi_\varphi\left(g \cdot(n\;k)\cdot
  h\right)\xi_\varphi,\xi_\varphi\right)=\varphi\left(g\cdot (n\;k)
  \cdot h\right)=\varphi\left((n\;m) \cdot g \cdot(n\;k)\cdot
  h\cdot (n\;m)\right)=\\
  \left(\pi_\varphi\left((n\;m) \cdot g (n\;k)\cdot
  h \cdot(n\;m)\right)\xi_\varphi,\xi_\varphi\right)=
  \left(\pi_\varphi\left(g\cdot (m\;k)\cdot
  h\right)\xi_\varphi,\xi_\varphi\right).
 \end{eqnarray*}
Approaching the limit as $m\to\infty$ we obtain the required
assertion.
\end{proof}

\begin{Lm}\label{abelian} The next relations hold true:
   \begin{itemize}
   \item[{\rm (1)}]\label{O1}$\mathcal{O}_k \mathcal{O}_n=\mathcal{O}_n \mathcal{O}_k$
   for all $
   k,n\in\mathbb{N}$;
   \item[{\rm (2)}]\label{O2}$\mathcal{O}_k
    \pi_\varphi\left(\gamma\right)=\pi_\varphi
   \left(\gamma\right) \mathcal{O}_k$
   for all $\gamma=\left(\gamma_1,\gamma_2,\ldots   \right)
   \in\Gamma_e^\infty$ such that  $\gamma_k=e $;
   \item[{\rm (3)}]\label{O3}$
   \pi_\varphi(s) \mathcal{O}_k=\mathcal{O}_{s(k)} \pi_\varphi(s)$
   for all  $s\in\mathfrak{S}_{\infty}$.
   \end{itemize}
\end{Lm}
 The proof follows immediately from definition $\mathcal{O}_k$
 (Proposition \ref{Prop19}). The details are left the reader.\qed

\medskip

 We will use the notation $\mathfrak{A}_j$ \label{definition mathfrakAj} for the
 $W^*-$algebra  generated by the operators
 $\pi_\varphi\left( \gamma\right)$, where
 $\gamma=\left(e,\cdots,e,\gamma_j,e,\cdots \right)$
 and $\mathcal{O}_j$. There is the natural isomorphism $\phi_{j,k}$ between
  $\mathfrak{A}_j$ and
  $\mathfrak{A}_k$ for any $k$ and $j$:
  \begin{eqnarray}\label{algebras}
   \phi_{j,k}:\mathfrak{A}_k\rightarrow
   \mathfrak{A}_j,\;\phi_{j,k}(a)=\pi_\varphi\left((k\;j)\right)
   a\pi_\varphi\left((k\;j)\right).
  \end{eqnarray}
  Observe that  $\left(\phi_{j,k}(a)\xi_\varphi,\xi_\varphi\right)=
  \left(a\xi_\varphi,\xi_\varphi\right)$ for all $k$, $j$ and $a\in \mathfrak{A}_k$.

  The next statement is the simple  technical generalization of proposition
\ref{multiplicativity}.
 \begin{Lm}\label{lemma22}
Let $s = \prod\limits_{p\in \mathbb{N}\diagup s}
 s_p $ be the decomposition of $s\in\mathfrak{S}_\infty$ into the
product  of cycles $s_p$, where $p\subset\mathbb{N}$ is the corresponding
orbit. Fix any finite collection $\left\{ U_j \right\}_{j=1}^N$ of the
elements from $\pi_\varphi \left( G \right)^{\prime\prime}$. If
$U_j\in\mathfrak{A}_j$ then
\begin{eqnarray}\label{gen mult formula}
  \left(\pi_\varphi(s)\prod \limits_j U_j\xi_\varphi,\xi_\varphi\right)=
  \prod\limits_{p\in \mathbb{N}/s}
  \left(\pi_\varphi(s_p)\prod \limits_{j\in p}
  U_j\xi_\varphi,\xi_\varphi\right).
 \end{eqnarray}
\end{Lm}
\begin{Prop}\label{fioncycle}
 Let $s_p\in\mathfrak{S}_\infty$ be the cyclic permutation on the set
 $p=\left\{k_1,k_2,\ldots,k_{|p\,|}\right\}\subset \mathbb{N},$
 where $k_l=s^{1-\,l}(k_1).$  If $U_{k_i}\in \mathfrak{A}_{k_i}$ for all
$k_i\in p$ then
   \begin{eqnarray}\label{mixture of U}
\begin{split}
     \left(\pi_\varphi(s_p)U_{k_1}U_{k_2}\cdots
     U_{k_{|p|}}\xi_\varphi,\xi_\varphi\right)\\=
\left(\phi_{k_{|p|}k_1}\left(U_{k_1}\right)\mathcal{O}_{k_{|p|}}
     \phi_{k_{|p|}k_2}\left(U_{k_2}\right)\mathcal{O}_{k_{|p|}}\cdots\mathcal{O}_{k_{|p|}}
     U_{k_{|p|}}\xi_\varphi,\xi_\varphi\right).
\end{split}
   \end{eqnarray}
\end{Prop}
 \begin{proof}
For convenience we suppose that $p=\{1,2,\ldots, n\}$ and
 \begin{eqnarray*}
  s_p(k)=\left\{\begin{array}{ll}
  k-1,&\textit{if }k>1\\n,&\textit{if }k=1
  \end{array}.\right.
\end{eqnarray*}
Since $s_p= (1\;n)(2\;n)\cdots (n-1\;\,n)$, we obtain
\begin{eqnarray*}
\begin{split}
&\left(\pi_\varphi \left( s_p \right) U_1U_2\cdots U_n\xi_\varphi,
 \xi_\varphi\right)\\
&=
\left(\pi_\varphi \left( (1\;n)(2\;n)\cdots(n-2\;n)  \right)
 U_1U_2\cdots\pi_\varphi \left((n-1\;\,n) \right)U_{n-1} U_n\xi_\varphi,
 \xi_\varphi\right)\\
 &=\left(\pi_\varphi \left( (1\;n)(2\;n)\cdots(n-2\;n)  \right)
 U_1U_2\cdots \phi_{n,n-1}\left(U_{n-1} \right)
 \pi_\varphi \left((n-1\;\,n)\right) U_n\xi_\varphi,
 \xi_\varphi\right).
\end{split}
\end{eqnarray*}
Hence, using $\mathfrak{S}_\infty$-invariance of $\varphi$ and lemma
\ref{abelian}, for any $N>n$ we have
\begin{eqnarray*}
\left(\pi_\varphi \left( s_p \right) U_1U_2\cdots U_n\xi_\varphi,
 \xi_\varphi\right)=\left(\pi_\varphi \left((n-1\;\,N) s_p (n-1\;\,N)\right) U_1U_2\cdots U_n\xi_\varphi,
 \xi_\varphi\right)\\
=\left(\pi_\varphi \left( (1\;n)(2\;n)\cdots(n-2\;n)  \right)
 U_1U_2\cdots \phi_{n,n-1}\left(U_{n-1} \right)
\pi_\varphi \left((N\;\,n)\right) U_n\xi_\varphi,
 \xi_\varphi\right).
\end{eqnarray*}
Approaching the limit as $N\to\infty$, we obtain
\begin{eqnarray*}
&\left(\pi_\varphi \left( s_p \right) U_1U_2\cdots U_n\xi_\varphi,
 \xi_\varphi\right)\\
&=
\left(\pi_\varphi \left( (1\;n)(2\;n)\cdots(n-2\;n)  \right)
 U_1U_2\cdots U_{n-2} \phi_{n,n-1}\left(U_{n-1} \right)
\mathcal{O}_n U_n\xi_\varphi,
 \xi_\varphi\right).
\end{eqnarray*}
Since $ \phi_{n,n-1}\left(U_{n-1} \right) \mathcal{O}_n$, then, by the
obvious induction, we have
\begin{eqnarray*}
&\left(\pi_\varphi \left( s_p \right) U_1U_2\cdots U_n\xi_\varphi,
 \xi_\varphi\right)\\
&=
\left(\phi_{n,1}\left(U_1 \right)
\mathcal{O}_n \phi_{n,2}\left(U_{2} \right)
\mathcal{O}_n \cdots  \phi_{n,n-2}\left(U_{n-2} \right)
\mathcal{O}_n \phi_{n,n-1}\left(U_{n-1} \right)
\mathcal{O}_n U_n\xi_\varphi,
 \xi_\varphi\right).
\end{eqnarray*}
\end{proof}

 The next statement is an analogue of Theorem 1 from
 {\cite{Ok2}}.
\begin{Lm}\label{discrete}
 Let $[a,b]$ belongs to $[-1,0]$ or
 $[0,1]$. with the property
. Denote by $E_{[a,b]}^{(i)}$ the spectral projection of
  self-adjoint operator $\mathcal{O}_i$.
  If  $\;\min\left\{ |a|,|b|  \right\}>\varepsilon>0$ then
$\left( E_{[a,b]}^{(i)}
\xi_\varphi,\xi_\varphi\right)^2\geq\varepsilon
 \left(E_{[a,b]}^{(i)}
 \xi_\varphi,\xi_\varphi\right)$.
 \end{Lm}
 This result may be proved in much the same way as
   theorem 1 from {\cite{Ok2}}. For convenience  we give below  the
  full  proof of  lemma \ref{discrete}.
 \begin{proof}
 Using Lemma {\ref{Fn}}, we
 have
 \begin{eqnarray}\label{est}
 \begin{split}
 &\left|\left( \pi_\varphi\left( (i,i+1)\right)
 E_{[a,b]}^{(i)}
 \xi_\varphi,
 E_{[a,b]}^{(i)}
 \xi_\varphi \right)\right|=\\
 &\left| \left( \mathcal{O}_i E_{[a,b]}^{(i)}\xi_\varphi,
  E_{[a,b]}^{(i)} \xi_\varphi \right)\right|\geqslant\varepsilon
 \left| \left( E_{[a,b]}^{(i)}\xi_\varphi, \xi_\varphi \right)\right|.&
 \end{split}
 \end{eqnarray}
 Hence, applying
 (\ref{algebras}) and lemma \ref{abelian}, we obtain
 \begin{eqnarray*}
&  E_{[a,b]}^{(i)}\pi_\varphi\left( (i,i+1)\right) E_{[a,b]}^{(i)}=
 E_{[a,b]}^{(i)} E_{[a,b]}^{(i+1)}
 \pi_\varphi\left( (i,i+1)\right)
 = \\&E_{[a,b]}^{(i)} E_{[a,b]}^{(i+1)}\pi_\varphi\left((i,i+1)\right)
 E_{[a,b]}^{(i)} E_{[a,b]}^{(i+1)}.
\end{eqnarray*}
Therefore,
 \begin{eqnarray*}
 &\left|\left( \pi_\varphi\left( (i,i+1)\right)
E_{[a,b]}^{(i)}
 \xi_\varphi,
 E_{[a,b]}^{(i)}
 \xi_\varphi \right)\right|\\
 &=\left|\left( \pi_\varphi\left( (i,i+1)\right)
E_{[a,b]}^{(i)} E_{[a,b]}^{(i+1)}
 \xi_\varphi, E_{[a,b]}^{(i)} E_{[a,b]}^{(i+1)}\xi_\varphi \right)\right|\\
 &\leq \left|\left(
E_{[a,b]}^{(i)} E_{[a,b]}^{(i+1)}
 \xi_\varphi, \xi_\varphi \right)\right|\stackrel{(Lemma
\ \ref{lemma22})}{=}
 \left(E_{[a,b]}^{(i)} E_{[a,b]}^{(i+1)}\xi_\varphi, \xi_\varphi \right)^2.
 \end{eqnarray*}
Hence, using ({\ref{est}}), we obtain the statement of  lemma
\ref{discrete}.
 \end{proof}

  Let
 $P_0^{(i)}$ be the orthogonal projection on
 ${\rm Ker}\,\mathcal{O}_i$.
  Put  $P_{\pm}^{(i)}=I-P_0^{(i)}$.

 \begin{Lm}\label{separating}
  Vector $\xi_\varphi$ is separating for  $w^*$-algebra
  $P^{(j)}_{\pm}\mathfrak{A}_j
  P_{\pm}^{(j)}$.
 \end{Lm}
 \begin{proof}
Let  $V\in P^{(j)}_{\pm}\mathfrak{A}_j
  P_{\pm}^{(j)}$ and let  $V\xi_\varphi=0$.
 It suffices to show that
\begin{eqnarray}\label{3.10aa}
\left(\pi_\varphi\left(g
\right)\xi_\varphi,\mathcal{O}_jV^*\pi_\varphi\left(h
\right)\xi_\varphi\right)=0 \;\text{ for all } \;g,h\in G.
\end{eqnarray}

First we note that, by $\mathfrak{S}_\infty$-invariance $\varphi$,
\begin{eqnarray}\label{inva}
\pi_\varphi\left(s \right) V\pi_\varphi\left(s^{-1}
\right)\xi_\varphi=0 \text{ for all } s\in\mathfrak{S}_\infty.
\end{eqnarray}
Further, if $g\in G_N$ then for all $n>N$
\begin{eqnarray*}
\pi_\varphi\left(\left(j\;n
\right)\right)V^*\pi_\varphi\left(\left(j\;n
\right)\right)\pi_\varphi\left(g\right)= \pi_\varphi\left(g\right)
 \pi_\varphi\left(\left(j\;n \right)\right)
V^*\pi_\varphi\left(\left(j\;n \right)\right).
\end{eqnarray*}
Hence, using definition of $\mathcal{O}_j$ (see proposition
\ref{Prop19}),
\begin{eqnarray*}
\begin{split}
\left(\pi_\varphi\left(g
\right)\xi_\varphi,\mathcal{O}_jV^*\pi_\varphi\left(h
\right)\xi_\varphi\right)=\lim\limits_{n\to\infty}
\left(\pi_\varphi\left(g
\right)\xi_\varphi,\pi_\varphi\left(\left(j\;n
\right)\right)V^*\pi_\varphi\left(h \right)\xi_\varphi\right)\\
=\lim\limits_{n\to\infty}\left(\pi_\varphi\left(\left(j\;n
\right)\right)V\pi_\varphi\left(\left(j\;n
\right)\right)\xi_\varphi,\pi_\varphi\left(g^{-1}
\right)\pi_\varphi\left(\left(j\;n \right)\right)\pi_\varphi\left(h
\right)\xi_\varphi\right)\stackrel{(\ref{inva})}{=}0.
\end{split}
\end{eqnarray*}
Thus (\ref{3.10aa}) is proved.
 \end{proof}
 The following statement is well known for the case of separating vector $\xi_\varphi$ (see
 {\cite{Ok2}}).
 In our case it follows from lemmas {\ref{discrete}} and {\ref{separating}}.
 \begin{Co}\label{co26}
 There exist at most countable set of numbers
 $\alpha_i$ from $[-1,0)\cup (0,1]$
 and a set of the pairwise orthogonal projections $
 \left\{P_{\alpha_i}^{(j)}\right\}
 \subset \mathfrak{A}_j$
 such that
 \begin{eqnarray}
 \mathcal{O}_j=P_0^{(j)}+\sum\limits_i \alpha_i P_{\alpha_i}^{(j)}.
\end{eqnarray}
  \end{Co}
\begin{Lm}\label{lemma27aa}
Let $\alpha,\beta\in {\rm Spectrum} \;\mathcal{O}_j$. If $\alpha \beta <0$ then
$P_{\alpha}^{(j)}\mathfrak{A}_jP_{\beta}^{(j)}=0$.
\end{Lm}
\begin{proof}
By lemma \ref{separating}, it suffices to show that
\begin{eqnarray}\label{103aa}
P_{\alpha}^{(j)}UP_{\beta}^{(j)}\xi_\varphi =0\;\text{ for all }
 U\in\mathfrak{A}_j.
\end{eqnarray}
First we note that
\begin{eqnarray*}
\left\|P_{\alpha}^{(j)}UP_{\beta}^{(j)}\xi_\varphi\right\|^2=
\left( P_{\beta}^{(j)}U^*P_{\alpha}^{(j)}UP_{\beta}^{(j)}\xi_\varphi,
\xi_\varphi\right)=\frac{1}{\alpha}\left( P_{\beta}^{(j)}U^*P_{\alpha}^{(j)}
\mathcal{O}_jUP_{\beta}^{(j)}\xi_\varphi,\xi_\varphi\right).
\end{eqnarray*}
Hence, using proposition \ref{fioncycle}, we receive
\begin{eqnarray}
\left\|P_{\alpha}^{(j)}UP_{\beta}^{(j)}\xi_\varphi\right\|^2=
\frac{1}{\alpha}\left( P_{\beta}^{(j)}U^*P_{\alpha}^{(j)}
\pi_\varphi \left(\left(j\;\,j+1 \right) \right)
P_{\alpha}^{(j)}UP_{\beta}^{(j)}\xi_\varphi,\xi_\varphi\right).
\end{eqnarray}
It follows from lemma \ref{abelian} that
\begin{eqnarray*}
\left\|P_{\alpha}^{(j)}UP_{\beta}^{(j)}\xi_\varphi\right\|^2=
\frac{1}{\alpha}\left( P_{\beta}^{(j)}U^*P_{\alpha}^{(j)}
\phi_{j+1,j}\left(P_{\alpha}^{(j)}UP_{\beta}^{(j)} \right)
\pi_\varphi \left(\left(j\;\,j+1 \right) \right)\xi_\varphi,
\xi_\varphi\right)\\
=\frac{1}{\alpha}\left(
\phi_{j+1,j}\left(P_{\alpha}^{(j)}UP_{\beta}^{(j)} \right)
P_{\beta}^{(j)}U^*P_{\alpha}^{(j)}
\pi_\varphi \left(\left(j\;\,j+1 \right) \right)\xi_\varphi,
\xi_\varphi\right)\\
=\frac{1}{\alpha}\left(
\phi_{j+1,j}\left(P_{\alpha}^{(j)}UP_{\beta}^{(j)} \right)
\pi_\varphi \left(\left(j\;\,j+1 \right) \right)
\phi_{j+1,j}\left(P_{\beta}^{(j)}U^*P_{\alpha}^{(j)}\right)
\xi_\varphi,\xi_\varphi\right)\\
=\frac{1}{\alpha}\left(
P_{\alpha}^{(j)}UP_{\beta}^{(j)}
\pi_\varphi \left(\left(j\;\,j+1 \right) \right)
P_{\beta}^{(j)}U^*P_{\alpha}^{(j)}
\xi_\varphi,\xi_\varphi\right)\\
\stackrel{\text{proposition \ref{fioncycle}}}{=}
\frac{1}{\alpha}\left(
P_{\alpha}^{(j)}UP_{\beta}^{(j)}
\mathcal{O}_j
P_{\beta}^{(j)}U^*P_{\alpha}^{(j)}
\xi_\varphi,\xi_\varphi\right)=
\frac{\beta}{\alpha}\left(
P_{\alpha}^{(j)}UP_{\beta}^{(j)}
U^*P_{\alpha}^{(j)}
\xi_\varphi,\xi_\varphi\right)\leq0.
\end{eqnarray*}
Therefore, (\ref{103aa}) holds true.
\end{proof}

The next assertion is an analogue
 of the theorem 2 from {\cite{Ok2}}.
 \begin{Lm}\label{P_alpha}
 Let $\alpha\neq 0$ be  the eigenvalue  of operator $\mathcal{O}_j$
 and let $P_\alpha^{(j)}$ be the corresponding spectral projection.
 Take any orthogonal projection
  $P\in P_\alpha^{(j)}\mathfrak{A}_jP_\alpha^{(j)}$ and put $\nu(P)=
 \left( P\xi_\varphi,\xi_\varphi\right)/|\alpha|$.
  Then
 $\nu(P)\in\mathbb{N}\cup\{0\}$.
 \end{Lm}
 \begin{proof}
 We use the arguments of
 Kerov, Olshanski, Vershik {\cite{KOV}} and Okounkov
 {\cite{Ok2}}. Let $j=1$.

 First consider the case $\alpha>0$.  For $n\in\mathbb{N}$
 put $\eta_n=\prod_{m=0}^{n-1} \phi_{1+m,1}(P)\xi_\varphi$.
 Let $s\in\mathfrak{S}_n$.
 In each orbit $p\in\mathbb{N}/s$ fix number
 $\mathfrak{s}(p)$.
  Since $\prod_{m=0}^{n-1} \phi_{1+m,1}(P)$
  is an orthogonal projection and
  \begin{eqnarray*}
  \pi_\varphi(s)\cdot\prod_{m=0}^{n-1} \phi_{1+m,1}(P)
  =\prod_{m=0}^{n-1} \phi_{1+m,1}(P)\cdot\pi_\varphi(s),
\end{eqnarray*}
   then we
  have
 \begin{eqnarray}\label{*}\begin{split}
 &\left(\pi_\varphi(s)\eta_n,\eta_n\right)=
 \left(\pi_\varphi(s)\prod_{m=0}^{n-1} \phi_{1+m,1}(P)\xi_\varphi,
\xi_\varphi\right)\\
& \stackrel{\text{lemma \ref{lemma22}}}{=}
\prod\limits_{p\in\left\{\mathbb{N}/s:p\,\subset [1,n]\right\}}
\left(\pi_\varphi(s_p)\prod\limits_{k\in p}
\phi_{k,j}(P)\xi_\varphi,
\xi_\varphi\right)\\
&\stackrel{\text{prop \ref{fioncycle}}}{=}
\prod\limits_{p\in\left\{\mathbb{N}/s:p\,\subset [1,n]\right\}}
\left(\phi_{{\mathfrak{s}(p)},1}\left(P\right)\cdot
\mathcal{O}_{\mathfrak{s}(p)}\cdot
\phi_{{\mathfrak{s}(p)},1}\left(P\right)\cdot\mathcal{O}_{\mathfrak{s}(p)}\cdots
\mathcal{O}_{\mathfrak{s}(p)}\cdot\phi_{{\mathfrak{s}(p)},1}\left(P\right)\xi_\varphi,\xi_\varphi\right)\\
 & =\prod\limits_{p\in\left\{\mathbb{N}/s:p\,\subset [1,n]\right\}}\alpha^{|p|-1}
\left(\phi_{{\mathfrak{s}(p)},1}\left(P\right)\xi_\varphi,\xi_\varphi\right)
 =\alpha^n\nu^{\,l(s)},
\end{split}
\end{eqnarray} where $l(s)$ is the number of cycles in the decomposition of permutation $s.$

Now define orthogonal projection
$Alt(n)\in\pi_\varphi\left(\mathfrak{S}_\infty
\right)^{\prime\prime} \subset\pi_\varphi(G)^{\prime\prime} $ by
\begin{eqnarray}
\label{Alt}
Alt(n)=\frac{1}{n!}\sum\limits_{s\in\mathfrak{S}_n}{\rm
sign}\,(s) \;\pi_\varphi(s).
\end{eqnarray}
 Using (\ref{*}), we obtain:
\begin{eqnarray}\label{Alt eta}
\left(Alt(n)\eta_n,\eta_n\right)=\alpha^n\sum\limits_{s\in\mathfrak{S}_n}
{\rm sign}\,(s)\;\nu^{\,l(s)}.
\end{eqnarray}
In the same way  as in \cite{Ok2}, applying
 equality:
 \begin{eqnarray*}\label{Alt equality}
 \sum\limits_{s\in\mathfrak{S}_n}{\rm
sign}\,(s)\;\nu^{\,l(s)}=\nu(\nu-1)\cdots(\nu-n+1),
 \end{eqnarray*}
we have
\begin{eqnarray}
0\leq
\left(\pi_\varphi(s)\eta_n,\eta_n\right)=\nu(\nu-1)\cdots(\nu-n+1).
\end{eqnarray}
  Therefore, $\nu\in \mathbb{N}\cup\{0\}$.

  The same proof remains for $\alpha<0$. In above reasoning
 operator $Alt(n)$ it is necessary to replace  by
$Sym(n)=\frac{1}{n!}\sum\limits_{s\in\mathfrak{S}_n}\pi_\varphi(s)$.
\end{proof}
For $\alpha\in {\rm Spectrum}\, \mathcal{O}_j$ denote by $P_\alpha^{(j)}$
the corresponding spectral projection (see corollary \ref{co26} ). It
follows from lemmas \ref{separating} and \ref{P_alpha} that for $\alpha
\neq 0$ $w^*$-algebra $P_\alpha^{(j)}\mathfrak{A}_jP_\alpha^{(j)}$ is
finite dimensional. Therefore, there exists finite collection $
\left\{P_{\alpha,i}^{(j)}\right\}_{i=1}^{n_\alpha }\subset P_\alpha^{(j)}
\mathfrak{A}_j P_\alpha^{(j)} $ of the {\it pairwise orthogonal}
projections with the properties:
\begin{eqnarray}\label{full system}
\begin{split}
&P_{\alpha,i}^{(j)}\xi_\varphi \neq 0\;\text{ and
}\;P_{\alpha,i}^{(j)} \;\text{ is minimal for all }
i=1,2,\ldots,\;n_\alpha;\\
&\sum\limits_{i=1}^{n_\alpha }P_{\alpha,i}^{(j)}=P_\alpha^{(j)}.
\end{split}
\end{eqnarray}
\begin{Prop}\label{P_zero}
Let $\mathcal{O}_j=P_0^{(j)}+\sum\limits_i \alpha_i
P_{\alpha_i}^{(j)}$. Put $P_+^{(j)}=\sum\limits_{i:\alpha_i>0}
P_{\alpha_i}^{(j)}$,  $P_-^{(j)}=\sum\limits_{i:\alpha_i<0}
P_{\alpha_i}^{(j)}$ and $P_\pm^{(j)}=P_+^{(j)}+P_-^{(j)}$.
 Then for each $U\in \mathfrak{A}_j$
\begin{eqnarray*}P_\pm^{(j)}UP_0^{(j)}\xi_\varphi=0.
\end{eqnarray*}
\end{Prop}
\begin{proof}
It is suffice to prove that $P_\alpha^{(j)}U P_0^{(j)}\xi_\varphi
=0$ for all nonzero $\alpha \in {\rm Spectrum}\,\mathcal{O}_j$. But
this fact follows from the next relations:
\begin{eqnarray*}
\left(P_\alpha^{(j)}U P_0^{(j)}\xi_\varphi,P_\alpha^{(j)}U
P_0^{(j)}\xi_\varphi \right)=\left(P_0^{(j)}U^* P_\alpha^{(j)}U
P_0^{(j)}\xi_\varphi,\xi_\varphi \right)\\
=\frac{1}{\alpha }\left(P_0^{(j)}U^*\mathcal{O}_j P_\alpha^{(j)}U
P_0^{(j)}\xi_\varphi,\xi_\varphi \right)\\
\stackrel{lemma\,\ref{Fn}}{=} \frac{1}{\alpha
}\left(P_0^{(j)}U^*\pi_\varphi \left(\left(j\;\,j+1 \right) \right)
P_\alpha^{(j)}U P_0^{(j)}\xi_\varphi,\xi_\varphi \right)\\
=\frac{1}{\alpha }\left(P_0^{(j)}U^*
P_\alpha^{(j+1)}\cdot\phi_{j+1,j}(U)\cdot P_0^{(j+1)}\pi_\varphi
\left(\left(j\;\,j+1 \right)\right) \xi_\varphi,\xi_\varphi
\right)\\
=\frac{1}{\alpha }\left(P_\alpha^{(j+1)}\cdot\phi_{j+1,j}(U)\cdot
P_0^{(j+1)}\cdot P_0^{(j)}U^* \pi_\varphi \left(\left(j\;\,j+1
\right)\right) \xi_\varphi,\xi_\varphi \right)\\
 =\frac{1}{\alpha }\left(P_\alpha^{(j+1)}\cdot\phi_{j+1,j}(U)\cdot
P_0^{(j+1)}\cdot  \pi_\varphi \left(\left(j\;\,j+1 \right)\right)
P_0^{(j+1)}\cdot \phi_{j+1,j}\left(U^*\right)\xi_\varphi,\xi_\varphi
\right)\\
\stackrel{lemma\,\ref{Fn}}{=}
 \frac{1}{\alpha }\left(P_\alpha^{(j+1)}\cdot\phi_{j+1,j}(U)\cdot
P_0^{(j+1)}\cdot\mathcal{O}_{j+1} \cdot P_0^{(j+1)}\cdot
\phi_{j+1,j}\left(U^*\right)\xi_\varphi,\xi_\varphi \right)=0.
\end{eqnarray*}
\end{proof}

Put $\mathbb{H}_{reg}^{(j)}= \left[ \mathfrak{A}_j
P_0^{(j)}\xi_\varphi \right]$ and $\mathbb{H}_{\pm}^{(j)}=
\left[\mathfrak{A}_j P_\pm^{(j)}\xi_\varphi  \right]$. The next
assertion follows from the previous proposition.
\begin{Co}\label{Co29}
\begin{itemize}
    \item [{\rm (a)}] Subspaces $\mathbb{H}_{reg}^{(j)}$ and $\mathbb{H}_{\pm}^{(j)}$ are
orthogonal for each $j\in \mathbb{N}$;
    \item [{\rm (b)}] if
     $\sum\limits_{\alpha\in {\rm Spectrum}\,\mathcal{O}_j:\alpha \neq 0}
    |\alpha| \cdot\nu \left(P_\alpha^{(j)}  \right)=1$ (see lemma \ref{P_alpha})
     then $P_0^{(j)}\xi_\varphi=0$.
\end{itemize}
\end{Co}
\begin{proof}
Property {\rm (a)} at once follows from proposition \ref{P_zero}. To
prove {\rm (b)} we note that $1=\left\|P_0^{(j)}\xi_\varphi
\right\|^2+\sum\limits_{\alpha\in {\rm
Spectrum}\,\mathcal{O}_j:\alpha \neq 0}
    \left\|P_\alpha^{(j)}\xi_\varphi\right\|^2 $ $\stackrel{\text{lemma \ref{P_alpha}}}{=}$
$\left\|P_0^{(j)}\xi_\varphi \right\|^2\\+\sum\limits_{\alpha\in {\rm
Spectrum}\,\mathcal{O}_j,\alpha \neq 0}\alpha \cdot\nu \left(P_\alpha^{(j)}
\right)$. Therefore, $\left\|P_0^{(j)}\xi_\varphi \right\|^2=0$.
\end{proof}
\begin{Lm}
$\left(U\mathcal{O}_jVP_0^{(j)}\xi_\varphi,P_0^{(j)}\xi_\varphi
\right)=0$ for all $U,V\in\mathfrak{A}_j$.
\end{Lm}

The proof follows from the next relations:
\begin{eqnarray*}
\left(U\mathcal{O}_jVP_0^{(j)}\xi_\varphi,P_0^{(j)}\xi_\varphi
\right)\stackrel{\text{lemma\,\ref{Fn}}}{=}
 \left(U\cdot\pi_\varphi \left(\left(j\;\,j+1 \right) \right)\cdot VP_0^{(j)}\xi_\varphi,P_0^{(j)}\xi_\varphi
\right)\\
=\left(P_0^{(j)}\cdot U\cdot \phi_{j+1,j}(V)\cdot
P_0^{(j+1)}\cdot\pi_\varphi \left(\left(j\;\,j+1 \right)
\right)\xi_\varphi,\xi_\varphi \right)\\
 =\left( \phi_{j+1,j}(V)\cdot
P_0^{(j+1)}\cdot P_0^{(j)}\cdot U\cdot\pi_\varphi
\left(\left(j\;\,j+1 \right) \right)\xi_\varphi,\xi_\varphi
\right)\\
 =\left( \phi_{j+1,j}(V)\cdot
P_0^{(j+1)}\cdot\pi_\varphi \left(\left(j\;\,j+1 \right)
\right)\cdot P_0^{(j+1)}\cdot \phi_{j+1,j}(U)\xi_\varphi,\xi_\varphi
\right)\\
\stackrel{\text{lemma\,\ref{Fn}}}{=} \left( \phi_{j+1,j}(V)\cdot
P_0^{(j+1)}\cdot\mathcal{O}_{j+1}\cdot P_0^{(j+1)}\cdot
\phi_{j+1,j}(U)\xi_\varphi,\xi_\varphi \right)=0.
\end{eqnarray*}\qed
\begin{Prop}\label{diff_Pa}
Let $ \left\{P_{\alpha,i}^{(j)} \right\}_{i=1}^{n_\alpha}$ $
\left(\alpha \in \left\{{\rm Spectrum}\,\mathcal{O}_j
\right\}\setminus 0\right)$ are the same as in (\ref{full system}).
If $P_{\alpha,i}^{(j)}\cdot P_{\beta,k}^{(j)}=0$ then $\left(
P_{\alpha,i}^{(j)}\cdot U\cdot
P_{\beta,k}^{(j)}\xi_\varphi,\xi_\varphi \right)=0$ for all
$U\in\mathfrak{A}_j$.
\end{Prop}
\begin{proof}
  The statement follows from the next relations:
  \begin{eqnarray*}
     &\left(P_{\alpha,i}^{(j)}\cdot U\cdot
P_{\beta,k}^{(j)}\xi_\varphi,\xi_\varphi \right)= \frac{1}{\alpha}
\left( P_{\alpha,i}^{(j)}\cdot\mathcal{O}_j\cdot U\cdot
P_{\beta,k}^{(j)}\xi_\varphi,\xi_\varphi \right)\\
         & \stackrel{\text{lemma\,\ref{Fn}}}{=}\frac{1}{\alpha}
\left( P_{\alpha,i}^{(j)}\cdot\pi_\varphi \left(\left(j\;\,j+1
\right) \right)\cdot U\cdot P_{\beta,k}^{(j)}\xi_\varphi,\xi_\varphi
\right)=\\
 &\frac{1}{\alpha}\left(P_{\alpha,i}^{(j)}\cdot \phi_{j+1,j}(U)\cdot
P_{\beta,k}^{(j+1)}\cdot\pi_\varphi \left(\left(j\;\,j+1 \right)
\right)\xi_\varphi,\xi_\varphi \right)\\
 &=\frac{1}{\alpha}
 \left(
 \phi_{j+1,j}(U)\cdot
P_{\beta,k}^{(j+1)}\cdot P_{\alpha,i}^{(j)}\cdot\pi_\varphi
\left(\left(j\;\,j+1 \right) \right)\xi_\varphi,\xi_\varphi
\right)\\
&=\frac{1}{\alpha} \left(  \phi_{j+1,j}(U)\cdot
P_{\beta,k}^{(j+1)}\cdot \pi_\varphi \left(\left( j\;\,j+1
\right)\right)
 \cdot P_{\alpha,i}^{(j+1)}\xi_\varphi,\xi_\varphi
\right)\\
     &\stackrel{lemma\,\ref{Fn}}{=}
     \frac{1}{\alpha}\left(  \phi_{j+1,j}(U)\cdot
P_{\beta,k}^{(j+1)}\cdot \mathcal{O}_{j+1}
 \cdot P_{\alpha,i}^{(j+1)}\xi_\varphi,\xi_\varphi
\right) \\
&=\left(  \phi_{j+1,j}(U)\cdot P_{\beta,k}^{(j+1)}\cdot
 P_{\alpha,i}^{(j+1)}\xi_\varphi,\xi_\varphi
\right)=0.
  \end{eqnarray*}
\end{proof}
Now we give important
\begin{Co} \label{Co32}
Let $P_+^{(j)}$ and  $P_-^{(j)}$ are the same as in proposition
\ref{P_zero}. Then subspaces $ \left[
\mathfrak{A}_jP_+^{(j)}\xi_\varphi \right]$ and $\left[
\mathfrak{A}_j P_-^{(j)}\xi_\varphi \right]$ are orthogonal.
\end{Co}

\begin{Prop}\label{Prop33}
Let $ \left\{P_{\alpha,i}^{(j)} \right\}_{i=1}^{n_\alpha}$ $
\left(\alpha \in \left\{{\rm Spectrum}\,\mathcal{O}_j
\right\}\setminus 0\right)$ are the same as in proposition
\ref{diff_Pa}. If there exists unitary $U\in\mathfrak{A}_j$ such
that $U\cdot P_{\alpha,i}^{(j)} \cdot U^*= P_{\beta,k}^{(j)} $ then
$\frac{\left(P_{\alpha,i}^{(j)}\xi_\varphi , \xi_\varphi
\right)}{|\alpha| }=\frac{\left(P_{\beta,k}^{(j)}\xi_\varphi ,
\xi_\varphi \right)}{|\beta| }$.
\end{Prop}
\begin{proof}
Let $\kappa_\alpha =\left(P_{\alpha,i}^{(j)}\xi_\varphi ,
\xi_\varphi \right)/|\alpha|$ and $\kappa_\beta
=\left(P_{\beta,k}^{(j)}\xi_\varphi , \xi_\varphi \right)/|\beta|$.
By lemma \ref{P_alpha}, $\kappa_\alpha, \kappa_\beta \in\mathbb{N}$.
Suppose for the convenience that $j=1$. For any $n\in\mathbb{N}$,
using (\ref{Alt}) and (\ref{Alt eta}), we obtain
\begin{eqnarray}\label{rel109}
\begin{split}
\left(Alt(n)\prod\limits_{m=1}^n\phi_{m,1}\left(P_{\alpha,i}^{(1)}
\right)\xi_\varphi,
\prod\limits_{m=1}^n\phi_{m,1}\left(P_{\alpha,i}^{(1)}
\right)\xi_\varphi
\right)=|\alpha|^n\prod\limits_{m=0}^{n-1}\left(\kappa_\alpha -m
\right);\\
 \left(Alt(n)\prod\limits_{m=1}^n\phi_{m,1}\left(P_{\beta,k}^{(1)}
\right)\xi_\varphi,
\prod\limits_{m=1}^n\phi_{m,1}\left(P_{\beta,k}^{(1)}
\right)\xi_\varphi
\right)=|\beta|^n\prod\limits_{m=0}^{n-1}\left(\kappa_\beta -m
\right).
\end{split}
\end{eqnarray}
This implies for $n=\kappa_\alpha+1$ that
\begin{eqnarray}\label{rel110}
\left(Alt(\kappa_\alpha+1)
\prod\limits_{m=1}^{\kappa_\alpha+1}\phi_{m,1}\left(P_{\alpha,i}^{(1)}
\right)\xi_\varphi,
\prod\limits_{m=1}^{\kappa_\alpha+1}\phi_{m,1}\left(P_{\alpha,i}^{(1)}
\right)\xi_\varphi \right)=0.
\end{eqnarray}
Further, applying relation
\begin{eqnarray*}
Alt(n)\cdot \prod\limits_{m=1}^n \phi_{m,1}(a)=\prod\limits_{m=1}^n
\phi_{m,1}(a)\cdot Alt(n) \;\;(\text{for all }a\in\mathfrak{A}_1),
\end{eqnarray*}
we get
\begin{eqnarray*}
&0\leq\left(Alt(\kappa_\alpha+1)\cdot
\prod\limits_{m=1}^{\kappa_\alpha+1}P_{\beta,k}^{(m)} \xi_\varphi,
\prod\limits_{m=1}^{\kappa_\alpha+1}P_{\beta,k}^{(m)}
\xi_\varphi \right)\\
& =\left(Alt(\kappa_\alpha+1)
\prod\limits_{m=1}^{\kappa_\alpha+1}P_{\alpha,i}^{(m)}
\phi_{m,1}\left(U^* \right)\xi_\varphi,
\prod\limits_{m=1}^{\kappa_\alpha+1}P_{\alpha,i}^{(m)}
\phi_{m,1}\left(U^* \right)\xi_\varphi \right)\\
 &=\left(Alt(\kappa_\alpha+1)
\prod\limits_{m=1}^{\kappa_\alpha+1}P_{\alpha,i}^{(m)}
\phi_{m,1}\left(U^* \right)\xi_\varphi,
\prod\limits_{m=1}^{\kappa_\alpha+1}\phi_{m,1}\left(U^*
\right)\xi_\varphi
\right)\\
&=\frac{1}{\alpha^{\kappa_\alpha+1} }
 \left(Alt(\kappa_\alpha+1)
\prod\limits_{m=1}^{\kappa_\alpha+1}P_{\alpha,i}^{(m)}\mathcal{O}_m
\phi_{m,1}\left(U^* \right)\xi_\varphi,
\prod\limits_{m=1}^{\kappa_\alpha+1}\phi_{m,1}\left(U^*
\right)\xi_\varphi
\right)\\
&\stackrel{\text{lemma \ref{Fn}}}{=}
\frac{1}{\alpha^{\kappa_\alpha+1} }
 \left(Alt(\kappa_\alpha+1)
\prod\limits_{m=1}^{\kappa_\alpha+1}P_{\alpha,i}^{(m)}\pi_\varphi
\left(\left(m\;\, \kappa_\alpha+1\right) \right)\phi_{m,1}\left(U^*
\right)\xi_\varphi,
\prod\limits_{m=1}^{\kappa_\alpha+1}\phi_{m,1}\left(U^*
\right)\xi_\varphi \right)\\
&=\frac{1}{\alpha^{\kappa_\alpha+1} }
 \left(Alt(\kappa_\alpha+1)
\prod\limits_{m=1}^{\kappa_\alpha+1}P_{\alpha,i}^{(m)}\phi_{m+\kappa_\alpha
+1,1}\left(U^* \right)\pi_\varphi \left(\left(m\;\,
\kappa_\alpha+1\right) \right)\xi_\varphi,
\prod\limits_{m=1}^{\kappa_\alpha+1}\phi_{m,1}\left(U^*
\right)\xi_\varphi \right)\\
&=\frac{1}{\alpha^{\kappa_\alpha+1} }
 \left(Alt(\kappa_\alpha+1)
\prod\limits_{m=1}^{\kappa_\alpha+1} P_{\alpha,i}^{(m)}\pi_\varphi
\left(\left(m\;\, \kappa_\alpha+1\right) \right)\xi_\varphi,
\prod\limits_{m=1}^{\kappa_\alpha+1}\phi_{m+\kappa_\alpha
+1,1}\left(U^* \right)\phi_{m,1}\left(U^* \right)\xi_\varphi
\right)\\
&\leq\frac{1}{|\alpha|^{\kappa_\alpha+1} } \left\|
Alt(\kappa_\alpha+1) \prod\limits_{m=1}^{\kappa_\alpha+1}
P_{\alpha,i}^{(m)}\pi_\varphi \left(\left(m\;\,
\kappa_\alpha+1\right) \right)\xi_\varphi\right\|\\
&=\frac{1}{|\alpha|^{\kappa_\alpha+1} } \left(Alt(\kappa_\alpha+1)
\prod\limits_{m=1}^{\kappa_\alpha+1} P_{\alpha,i}^{(m)}\pi_\varphi
\left(\left(m\;\, \kappa_\alpha+1\right)
\right)\xi_\varphi,\pi_\varphi \left(\left(m\;\,
\kappa_\alpha+1\right) \right)\xi_\varphi \right)^{1/2}\\
&\stackrel{\mathfrak{S}_\infty\text{-centrality of } \varphi }{=}
\frac{1}{|\alpha|^{\kappa_\alpha+1} } \left(Alt(\kappa_\alpha+1)
\prod\limits_{m=1}^{\kappa_\alpha+1}
P_{\alpha,i}^{(m)}\xi_\varphi,\xi_\varphi
\right)^{1/2}\stackrel{(\ref{rel110})}{=}0.
\end{eqnarray*}
Hence, applying (\ref{rel109}), we have
\begin{eqnarray*}
\left(Alt(\kappa_\alpha+1)\cdot
\prod\limits_{m=1}^{\kappa_\alpha+1}P_{\beta,k}^{(m)} \xi_\varphi,
\prod\limits_{m=1}^{\kappa_\alpha+1}P_{\beta,k}^{(m)} \xi_\varphi
\right)\\
 =|\beta|^{\kappa_\alpha +1 }\kappa_\beta \left(\kappa_\beta-1 \right)
 \left(\kappa_\beta-2\right)\left(\kappa_\beta-\kappa_a\right)=0.
\end{eqnarray*}
Therefore, $\kappa_\alpha \geq\kappa_\beta$. Similarly,
$\kappa_\alpha \leq\kappa_\beta$.
\end{proof}
\paragraph{The proof of theorem \ref{main}.}
Now we will give the description  of parameters
  $\left(A,\rho\right)$  from paragraph
\ref{paragraph2.1}, corresponding to $\varphi$.

First we describe  the structure of $w^*$-algebra $\widetilde{P}_\pm^{(j)}
\mathfrak{A}_j$, \footnote{ see page \pageref{definition mathfrakAj}  for
the definition of $\mathfrak{A}_j$}  where $\widetilde{P}_\pm^{(j)}$ is the
orthogonal projection of $\left[\mathfrak{A}_j\xi_\varphi \right]$ onto $
\left[ \mathfrak{A}_jP_\pm^{(j)}\xi_\varphi\right]$ \footnote{$P_\pm^{(j)}$
is defined in proposition \ref{P_zero}}.

Let $\mathcal{C}_\pm^{(j)}$ be the center of
$\widetilde{P}_\pm^{(j)} \mathfrak{A}_j$. Denote by
$c(P)\in\mathcal{C}_\pm^{(j)}$ the central support of projection
$P\in\widetilde{P}_\pm^{(j)} \mathfrak{A}_j$. Let us prove that
\begin{eqnarray}\label{central111}
c \left( P_\pm^{(j)}\right)=\widetilde{P}_\pm^{(j)}.
\end{eqnarray}
Indeed, if $F=\widetilde{P}_\pm^{(j)}-c\left(P_\pm^{(j)}\right)$,
then for all $B\in\mathfrak{A}_j$ we have
 $F BP_\pm^{(j)}\xi_\varphi$ $ =
 BFP_\pm^{(j)}\xi_\varphi=0$. Therefore, $F=0$.

 Since for any nonzero $\alpha \in \left\{{\rm Spectrum}\,\mathcal{O}_j\right\}\setminus 0$
  in $P_\alpha ^{(j)} \mathfrak{A}_jP_\alpha ^{(j)}$ there exists
  finite collection $
\left\{P_{\alpha,i}^{(j)}\right\}_{i=1}^{n_\alpha }$ of the {\it
minimal } projections with properties (\ref{full system}), then
$w^*$-algebra $P_\pm^{(j)}\mathfrak{A}_jP_\pm^{(j)}$ is
$*$-isomorphic to
 the direct sum of full matrix algebras.
 Thus, using (\ref{central111}), we  find the collection $ \left\{F_m \right\}_{m=1}^N$ of pairwise
 orthogonal projections from $\mathcal{C}_\pm^{(j)}$ such that
 $F_m\cdot\widetilde{P}_\pm^{(j)} \mathfrak{A}_j\cdot F_m$ is a factor
 of the type ${\rm I}_{k_m}$. Denote
 $F_m\cdot\widetilde{P}_\pm^{(j)} \mathfrak{A}_j\cdot
 F_m$ by $\mathcal{M}_{k_m}$. That is $P_\pm^{(j)}\mathfrak{A}_jP_\pm^{(j)}$ is
isomorphic to $\mathcal{M}_{k_1}\oplus \mathcal{M}_{k_2}\oplus\ldots$.
  Let
 $ \left\{e^{(m)}_{pq} \right\}_{p,q=1}^{k_m}$ be the matrix unit of
$\mathcal{M}_{k_m}$.
Without loss   of the generality we suppose that for certain $l_m\leq k_m$
\begin{eqnarray}\label{111}
\begin{split}
\bigcup\limits_m\left\{e^{(m)}_{pp}
\right\}_{p=1}^{l_m}\subset\bigcup\limits_{\alpha \in{\rm
Spectrum}\,\mathcal{O}_j,\,\alpha \neq
0}\left\{P_{\alpha,i}^{(j)}\right\}_{i=1}^{n_\alpha }\;\;\text{ and }\\
 \left\{\bigcup\limits_m\left\{e^{(m)}_{pp}
\right\}_{p=l_m+1}^{k_m}\right\}\, \bigcap \, \left\{\bigcup\limits_{\alpha \in{\rm
Spectrum}\,\mathcal{O}_j,\,\alpha \neq
0}\left\{P_{\alpha,i}^{(j)}\right\}_{i=1}^{n_\alpha }\right\}=\emptyset.
\end{split}
\end{eqnarray}
By lemmas \ref{separating}, \ref{P_alpha} and propositions \ref{diff_Pa},
\ref{Prop33}, minimal projections $\bigcup\limits_m\left\{e^{(m)}_{pp}
\right\}_{p=1}^{l_m}$ satisfy the next conditions
\begin{itemize}
  \item {\rm (a)} if $e^{(m)}_{pp}\cdot\mathcal{O}_j=\alpha_p\cdot
      e^{(m)}_{pp}$, where $\alpha_p\in{\rm Spectrum}\setminus 0$, then
      there exists natural $q_m$ such that
      $\frac{\left(e^{(m)}_{pp}\xi_\varphi,\xi_\varphi\right)}{\left|\alpha_p\right|
      }=q_m$ for all $p=1,2,\ldots, l_m$;
  \item {\rm (b)}  if $p\neq q $ then $\left(
      e^{(m)}_{pq}\xi_\varphi,\xi_\varphi\right)=0$ for all
      $p,q=1,2,\ldots, l_m$; $m=1,2,\ldots, N$.
\end{itemize}
Further, using (\ref{111}), for $p>l_m$ we have
\begin{eqnarray*}
e_{pp}^{(m)}\cdot P_0^{(j)}=e_{pp}^{(m)}.
\end{eqnarray*}
It follows from this and proposition \ref{P_zero} that
\begin{eqnarray}\label{0.112}
\left(e_{pq}^{(m)} \xi_\varphi,\xi_\varphi\right)=0\;\text{ for }\;\;
p=1,2,\ldots,l_m; \;\; q=l_m+1, l_m+2, \ldots,k_m.
\end{eqnarray}
Let us prove that
\begin{eqnarray}\label{0.113}
\left(e_{pq}^{(m)} \xi_\varphi,\xi_\varphi\right)=0\;\text{ for }\;\;
 p,q=l_m+1,l_m+2,\ldots,k_m.
\end{eqnarray}
For this it suffices to prove the next equality:
\begin{eqnarray}\label{0.114}
\left(e_{pp}^{(m)} \xi_\varphi,\xi_\varphi\right)=0\;\text{ for }\;\;
 p,q=l_m+1,l_m+2,\ldots,k_m.
\end{eqnarray}
Fix $p>l_m$. Applying proposition \ref{fioncycle}, we have
\begin{eqnarray*}
&\left(e_{pp}^{(m)} \xi_\varphi,\xi_\varphi\right)=\frac{1}{\alpha_1 }
  \left(e_{p1}^{(m)} \cdot\mathcal{O}_j\cdot e_{1p}^{(m)}\xi_\varphi,
  \xi_\varphi \right)\\
  &\stackrel{\text{proposition \ref{fioncycle}}}{=}
  \frac{1}{\alpha_1 }\left(\pi_\varphi \left(\left(j\,\;j+1 \right) \right)\cdot
  \phi_{j+1,j}\left(e_{p1}^{(m)}\right)\cdot e_{1p}^{(m)}\xi_\varphi,
  \xi_\varphi \right)\\
 & =\frac{1}{\alpha_1 }
  \left(\pi_\varphi \left(\left(j\,\;j+1 \right) \right)\cdot e_{1p}^{(m)}
  \cdot\phi_{j+1,j}\left(e_{p1}^{(m)}\right)\xi_\varphi,
  \xi_\varphi \right)\\
  &=\frac{1}{\alpha_1 }
  \left(
  \phi_{j+1,j}\left(e_{1p}^{(m)}\right)\cdot\pi_\varphi \left(\left(j\,\;j+1 \right) \right)
  \cdot\phi_{j+1,j}\left(e_{p1}^{(m)}\right)\xi_\varphi,
  \xi_\varphi \right)\\
  &=\frac{1}{\alpha_1 }
  \left(\pi_\varphi \left(\left(j\;\,n \right) \right)\cdot
  \phi_{j+1,j}\left(e_{1p}^{(m)}\right)\cdot\pi_\varphi \left(\left(j\,\;j+1 \right) \right)
  \cdot\phi_{j+1,j}\left(e_{p1}^{(m)}\right)\cdot
  \pi_\varphi \left(\left(j\;\,n \right) \right)\xi_\varphi,
  \xi_\varphi \right)\\
  &=\frac{1}{\alpha_1 }
  \left(
  \phi_{j+1,j}\left(e_{1p}^{(m)}\right)\cdot\pi_\varphi \left(\left(j+1\,\;n \right) \right)
  \cdot\phi_{j+1,j}\left(e_{p1}^{(m)}\right)\xi_\varphi,
  \xi_\varphi \right)\\
  &=\lim\limits_{n\to\infty}\frac{1}{\alpha_1 }
  \left(
  \phi_{j+1,j}\left(e_{1p}^{(m)}\right)\cdot\pi_\varphi \left(\left(j+1\,\;n \right) \right)
  \cdot\phi_{j+1,j}\left(e_{p1}^{(m)}\right)\xi_\varphi,
  \xi_\varphi \right)\\
  &=\frac{1}{\alpha_1 }
  \left(
  \phi_{j+1,j}\left(e_{1p}^{(m)}\right)\cdot\mathcal{O}_{j+1}
  \cdot\phi_{j+1,j}\left(e_{p1}^{(m)}\right)\xi_\varphi,
  \xi_\varphi \right)\\
  &=\frac{1}{\alpha_1 }\left(e_{1p}^{(m)}\cdot
  \mathcal{O}_j\cdot e_{p1}^{(m)}\xi_\varphi,\xi_\varphi\right)\stackrel{(\ref{111})}{=}0.
  \end{eqnarray*}
  Thus (\ref{0.114}) and (\ref{0.113}) are proved.

    Define $\widehat{\varphi }\in \pi_\varphi \left(G\right)^{\prime\prime}_*$
by $\widehat{\varphi}(a)=\left(a\xi_\varphi,\xi_\varphi\right)$. Denote by
$\mathbb{M}_{q_m}$ the algebra of all complex matrices and put $\mathcal{N}_m=
 \mathcal{M}_{k_m}\otimes \mathbb{M}_{q_m}$,
  $A^{(m)}=\sum\limits_{p=1}^{l_m}\alpha_p\cdot e_{pp}^{(m)}\in\mathcal{M}_{k_m}$
  (see property {\rm (a)} and (\ref{111})). Consider $w^*$-algebra
$\widetilde{\mathfrak{A}}_j=\left(\bigoplus\limits_{m=1}^N
 F_m\mathfrak{A}_jF_m\otimes\mathbb{M}_{q_m}\right)\bigoplus
 \left(I-\widetilde{P}_\pm^{(j)}\right)\mathfrak{A}_j$.
Observe that there exists the natural embedding
\begin{eqnarray}
\mathfrak{A}_j\ni a\stackrel{\mathfrak{i}}{\mapsto}
\sum\limits_{m=1}^N \left(F_maF_m\otimes I \right)+
 \left(I-\widetilde{P}_\pm^{(j)}\right)a\in \widetilde{\mathfrak{A}}_j.
\end{eqnarray}
   Now, using properties
   {\rm (a)}-{\rm (b)}, (\ref{0.112}) and
   (\ref{0.113}), we have for all
  $a\in  \mathfrak{A}_j$
  \begin{eqnarray}\label{0.116}
\widehat{\varphi }
\left(a \right)=\sum_{m=1}^N {\rm Tr}_m\left(a\left|A^{(m)} \right|\otimes I \right)+
\left(a\left(I-\widetilde{P}_\pm^{(j)} \right)\xi_\varphi,\xi_\varphi\right),
  \end{eqnarray}
  where ${\rm Tr}_m$ is ordinary trace\footnote{If $e$ is minimal
  projection from $\mathcal{N}_m$ then ${\rm Tr}_m(e)=1$. }
  on $\mathcal{N}_m$.

  Now we define parameters $\left\{\mathcal{H},A,\rho,\hat{\xi}\right\}$ from
  paragraph \ref{paragraph2.1} such that
  \begin{eqnarray}\label{0.117}
\varphi =\psi_A^\rho\;\;(\text{ see proposition \ref{Prop11a}}).
  \end{eqnarray}
  For this purpose we fix in each $\mathcal{N}_m=
 \mathcal{M}_{k_m}\otimes \mathbb{M}_{q_m}$ minimal projection $e_m$. Define
 state $f$ on $\widetilde{\mathfrak{A}}_j$ by
 \begin{eqnarray}
f\left(\widetilde{a}\right)=\sum\limits_{m=1}^N{\rm Tr}_m\left(
e_m\widetilde{a} e_m\right) \;\;\;
 \left(\widetilde{a}\in\widetilde{\mathfrak{A}}_j \right).
 \end{eqnarray}
 Let $\left(R_f,\mathcal{H}_f,\xi_f\right)$ be the corresponding GNS-representation of
 $\widetilde{\mathfrak{A}}_j$. Now we define $\mathcal{H}$ by
 \begin{eqnarray}
\mathcal{H}=\mathcal{H}_f\oplus
\left[\left(I-\widetilde{P}_\pm^{(1)}
\right)\mathfrak{A}_1\xi_\varphi \right]
\oplus
\left[\left(I-\widetilde{P}_\pm^{(2)}\right)\mathfrak{A}_2\xi_\varphi
\right]\oplus\ldots.
\end{eqnarray}
Representation $\rho$ acts on $\eta_p\in
\left[\left(I-\widetilde{P}_\pm^{(p)}\right)\mathfrak{A}_p\xi_\varphi
\right]$ as follows
\begin{eqnarray}
\rho\left(\gamma\right)\eta_p= \pi_\varphi
\left(\left(e,\ldots,\stackrel{p-th}{\gamma},e,\ldots  \right) \right)\eta_p.
\end{eqnarray}
If $\eta\in\mathcal{H}_f$ then
\begin{eqnarray}
\rho\left(\gamma\right)\eta=R_f\circ\mathfrak{i}\left(
 \pi_\varphi
\left(\left(e,\ldots,\stackrel{j-th}{\gamma},e,\ldots  \right)
\right)\right)\eta.
\end{eqnarray}
Operator $A$ is defined by
\begin{eqnarray}
A\eta=\left\{\begin{array}{ll}
R_f\circ\mathfrak{i}\left(\sum\limits_{m=1}^NA^{(m)} \right)\eta,
&\text{ if }\;\eta\in\mathcal{H}_f,\\
0,&\text{ if }\;\eta\in\left[\left(I-\widetilde{P}_\pm^{(p)}\right)\mathfrak{A}_p\xi_\varphi
\right].\end{array}\right.
\end{eqnarray}
In the case $\sum\limits_{\alpha \in{\rm
Spectrum}\,\mathcal{O}_j,\,\alpha \neq
0}|\alpha|\nu\left(P_\alpha^{(j)}\right)=\sum\limits_{m=1}^{N}
\sum\limits_{p=1}^{k_m}\left|\alpha_p\right|<1$
(see corollary \ref{Co29} and property {\rm(a)}) vector $\hat{\xi}$ is defined
 by
 \begin{eqnarray}
\hat{\xi}=\frac{\left(I-\widetilde{P}_\pm^{(1)}\right)\xi_\varphi }{
 \left\| \left(I-\widetilde{P}_\pm^{(1)}\right)\xi_\varphi \right\|}.
 \end{eqnarray}
 Now it follows from (\ref{0.116}) that for $a\in\mathfrak{A}_j$
 \begin{eqnarray}
 \begin{split}
&\widehat{\varphi}(a) ={\rm Tr}\left(R_f\left(\mathfrak{i}(a)\right)\cdot
\left|A\right| \right)\\
&+\left\| \left(I-\widetilde{P}_\pm^{(1)}\right)\xi_\varphi \right\|
\left(\pi_\varphi\left(\left(1\;\,j \right) \right)\cdot a\cdot\pi_\varphi\left(\left(1\;\,j \right) \right)\hat{\xi},\hat{\xi}\right).
\end{split}
 \end{eqnarray}
 Hence, applying lemma \ref{lemma22}, proposition \ref{fioncycle} and
  definition of $\psi_A^\rho$, we can to receive equality (\ref{0.117}).
 In particular, lemma \ref{lemma27aa} implies property {\rm(3)}
  from paragraph \ref{paragraph2.1}. \qed


{}

\medskip
{\it A.V. Dudko,  Department of Mathematics, University of Toronto,

 e-mail:
 artem.dudko@utoronto.ca

\medskip
N.I. Nessonov, B.Verkin Institute for Low Temperature Physics and Engineering
of the National Academy of Sciences of Ukraine,

e-mail: nessonov@ilt.kharkov.ua}
\end{document}